\documentclass{article}
\usepackage{amssymb}
\usepackage{amsmath,color}
\usepackage[hidelinks]{hyperref}

\usepackage[square,numbers]{natbib}         
\usepackage{ntheorem}
\newtheorem{lem}{Lemma}[section]
\newtheorem{cor}[lem]{Corollary}
\newtheorem{teo}[lem]{Theorem}
\newtheorem*{teoNon}{Theorem}
\newtheorem{os}[lem]{Remark}
\newtheorem{defi}[lem]{Definition}
\newtheorem{prop}[lem]{Proposition}

\newtheorem{assumption}[lem]{Assumption}

\newcommand{\qed}{\thinspace\null\nobreak\hfill\hbox{\vbox{\kern-.2pt\hrule
 height.2pt depth.2pt\kern-.2pt\kern-.2pt \hbox to2.5mm{\kern-.2pt\vrule
 width.4pt \kern-.2pt\raise2.5mm\vbox to.2pt{}\lower0pt\vtop
 to.2pt{}\hfil\kern-.2pt \vrule
 width.4pt \kern-.2pt}\kern-.2pt\kern-.2pt\hrule height.2pt depth.2pt
 \kern-.2pt}}\par\medbreak}

\newcommand{\R}{\mathbb{R}}

\newcommand{\C}{\mathbb{C}}

\newcommand{\Rp}{\textrm{\emph{Re}\,}}

\newcommand{\ds}{\displaystyle}

\date{}
\begin{document}

\title{Singular parabolic operators in the half-space with  boundary degeneracy: Dirichlet  and oblique derivative boundary conditions}
\author{L. Negro \thanks{Dipartimento di Matematica e Fisica  ``Ennio De
Giorgi'', Universit\`a del Salento, C.P.193, 73100, Lecce, Italy. email: luigi.negro@unisalento.it}}

\maketitle
\begin{abstract}
\noindent 
We study elliptic and parabolic problems governed by the  singular elliptic   operators 
\begin{align*}
	\mathcal L=y^{\alpha_1}\mbox{Tr }\left(QD^2_x\right)+2y^{\frac{\alpha_1+\alpha_2}{2}}q\cdot \nabla_xD_y+\gamma y^{\alpha_2} D_{yy}+y^{\frac{\alpha_1+\alpha_2}{2}-1}\left(d,\nabla_x\right)+cy^{\alpha_2-1}D_y-by^{\alpha_2-2}
\end{align*}
in the half-space $\R^{N+1}_+=\{(x,y): x \in \R^N, y>0\}$, under Dirichlet or oblique derivative boundary conditions. In the special case $\alpha_1=\alpha_2=\alpha$ the operator $\mathcal L$ takes the form
\begin{align*}
	\mathcal L&=y^{\alpha}\mbox{Tr }\left(AD^2\right)+y^{\alpha-1}\left(v,\nabla\right)-by^{\alpha-2},
\end{align*}
where $v=(d,c)\in\R^{N+1}$, $b\in\R$ and 
$
A=\left(
\begin{array}{c|c}
	Q  & { q}^t \\[1ex] \hline
	q& \gamma
\end{array}\right)$ is an elliptic matrix.
 We prove elliptic and
parabolic $L^p$-estimates and solvability for the associated problems. In the language of semigroup
theory, we prove that $\mathcal L$ generates an analytic semigroup, characterize its domain as a weighted
Sobolev space and show that it has maximal regularity.

\bigskip\noindent
Mathematics subject classification (2020): 35K67, 35B45, 47D07, 35J70, 35J75.
\par

\noindent Keywords: degenerate elliptic operators, boundary degeneracy, vector-valued harmonic analysis,  maximal regularity.
\end{abstract}

\section{Introduction}
In this paper we study solvability and regularity of elliptic and parabolic problems associated to the  degenerate   operators
\begin{align}\label{defL}
	\mathcal L=y^{\alpha_1}\mbox{Tr }\left(QD^2_x\right)+2y^{\frac{\alpha_1+\alpha_2}{2}}q\cdot \nabla_xD_y+\gamma y^{\alpha_2} D_{yy}+y^{\frac{\alpha_1+\alpha_2}{2}-1}\left(d,\nabla_x\right)+cy^{\alpha_2-1}D_y-by^{\alpha_2-2}
\end{align}
and $D_t- \mathcal L$ in the half-space $\R^{N+1}_+=\{(x,y): x \in \R^N, y>0\}$ or  in $(0, \infty) \times \R^{N+1}_+$ and under Dirichlet or  oblique derivative boundary conditions at $y=0$.

Here $v=(d,c)\in\R^{N+1}$ with $d=0$ if $c =0$, $b\in\R$ and 
$
A=\left(
\begin{array}{c|c}
	Q  & { q}^t \\[1ex] \hline
	q& \gamma
\end{array}\right)$ is 
a constant real  elliptic matrix. The real numbers $\alpha_1, \alpha_2$ satisfy $\alpha_2<2$ and $\alpha_2-\alpha_1<2$ but are not assumed to be nonnegative.  In the special case $\alpha_1=\alpha_2=\alpha$ the operator $\mathcal L$ takes the form
\begin{align}\label{defL alpha eq}
	\mathcal L&=y^{\alpha}\mbox{Tr }\left(AD^2\right)+y^{\alpha-1}\left(v,\nabla\right)-by^{\alpha-2}
\end{align}
whose coefficients are singular for $\alpha<0$ and degenerate for $\alpha>0$ at $y=0$.
This paper is the companion of \cite{MNS-Degenerate-Half-Space}, in which the same type of operators are considered but with $d=0$, $b=0$, and with Neumann boundary condition.

We write
$
B_y
$  to denote the 1-dimensional  Bessel operator $D_{yy}+\frac{c}{y}D_y$ and $L_y=D_{yy}+\frac{c}{y}D_y-\frac{b}{y^2}$;  note that $B_y$  is nothing but $L_y$ when $b=0$. With this notation  the   special cases where
\begin{equation} \label{Lpart}
\mathcal L=y^{\alpha_1}\Delta_x+y^{\alpha_2}B_y,\qquad \mathcal L=y^{\alpha_1}\Delta_x+y^{\alpha_2}L_y
\end{equation}
has been already studied in \cite{MNS-Caffarelli, MNS-CompleteDegenerate}.  The main novelty here consists in the presence of the mixed derivatives $2y^{\frac{\alpha_1+\alpha_2}{2}}q\cdot \nabla_xD_y$ and of the $x$-derivative  $y^{\frac{\alpha_1+\alpha_2}{2}-1}\left(d,\nabla_x\right)$ in the operator $\mathcal L$ which is a crucial step for treating degenerate operators in domains, through a localization procedure.  Surprisingly enough, the case $\alpha_1=\alpha_2$ implies all other cases by the change of variables described in Section \ref{Section Degenerate}. However this modifies the underlying measure and the procedure works if one is able to deal with the complete scale of $L^p_m$ spaces, where $L^p_m=L^p(\R^{N+1}_+; y^m dxdy)$.\\

The interest in this class of singular  operators has grown in the last decade as they appear extensively  in the literature in 
both pure and applied problems. The operators in \eqref{Lpart} are strongly connected with  nonlocal operators as they   play a major role in the investigation of the fractional powers of the Laplacian and of the Heat operator
through the  ``extension procedure" of Caffarelli and Silvestre, see \cite{Caffarelli-Silvestre} and \cite{Stinga-Torrea-Extension,GaleMiana-Extension,Arendt-ExtensionProblem} for a more general setting.
We refer also the reader to the introductions of \cite{dong2020parabolic,Vita2024schauder}  for some references
to related problems in probability, mathematical
finance and biology, porous media equations and in degenerate viscous Hamilton-Jacobi equations.  This type of singular  operators are also connected to the theory of geometric PDEs with edge singularities  \cite{MazzeoEdgeTheory} and  to the analysis of the regularity of the ratio of solutions to elliptic PDEs \cite{Vita2024Harnack}.\\

Our main results are Theorems \ref{Teo alpha Oblique}, \ref{complete},  \ref{complete-oblique} and \ref{complete dirichlet alpha diff}, where, in the language of semigroup theory, we prove that $\mathcal L$ generates an analytic semigroup on $L^p_m$, characterize its domain as a weighted Sobolev space and show that it has maximal parabolic regularity.  For reader's convenience we collect, in  Section \ref{Section main results},  the  main hypotheses  we assume  and also the main  results of the paper in the case $\alpha_1=\alpha_2$ as in  \eqref{defL alpha eq}, referring to  Section \ref{Consequenses} for their extension to  general $\alpha_1$, $\alpha_2$ as in \eqref{defL}.\\

 We prove both elliptic and parabolic estimates which, in the case $\alpha_1=\alpha_2=\alpha$, $b=0$ and oblique derivative boundary condition, read as
\begin{align} \label{closedness}
	\|y^{\alpha}D^2 u\|_{L^p_m}
	+\| y^{\alpha-1}v\cdot \nabla  u\|_{L^p_m} \le C\| \mathcal L u\|_{L^p_m},
\end{align}
and
\begin{equation} \label{closedness1}
	\|D_t u\|_{L^p_m}+\|\mathcal Lu\|_{L^p_m} \le C\| (D_t-\mathcal L) u\|_{L^p_m}, \quad 
\end{equation}
where the $L^p$ norms are taken over $\R_+^{N+1}$ and on $(0, \infty) \times \R_+^{N+1}$ respectively. Both the elliptic and parabolic estimates above  share the name ``maximal regularity" even though this term is often restricted to the parabolic case. Throughout the paper we keep this convention referring, in the statements of our results, to maximal regularity to denote the validity of the  parabolic estimates \eqref{closedness1} while the elliptic bounds as in  \eqref{closedness} will be characterized through the precise description of the domain of $\mathcal L$.   \\

Let us explain the meaning of the restrictions $\alpha_2<2$, $\alpha_2-\alpha_1<2$
considering first the case where $\alpha_1=\alpha_2=\alpha$, so that the unique requirement is $\alpha<2$. 

It turns out that when $\alpha \geq 2$ the problem is easily treated in the strip $\R^N\times [0,1]$ in the case of the Lebesgue measure, see \cite{FornMetPallScn5}, and all problems are due to the strong diffusion at infinity. The case $\alpha \geq 2$ in the strip $\R^N \times [1, \infty[$ requires therefore new investigation even though the 1-dimensional case is easily treated by the change of variables of Section \ref{Section Degenerate}.

When $\alpha_1 \neq \alpha_2$, the further restriction $\alpha_2-\alpha_1<2$ comes from the change of variables of Section \ref{Section Degenerate}, see Section \ref{Consequenses}.\\

Let us briefly describe the previous literature on these  operators. In \cite{MNS-Caffarelli, MNS-Caf-Schauder}  we considered the simplest case of  $\Delta_x+B_y$ making extensive use of  the commutative structure of the operator. The non-commutative case of $y^{\alpha_1}\Delta_x+y^{\alpha_2}B_y$ has been later faced in \cite{MNS-CompleteDegenerate}.  
Another source of non-commutativity comes from the presence of mixed derivatives. In \cite{MNS-Degenerate-Half-Space, MNS-Singular-Half-Space, Negro-Spina-SingularKernel} we treated the operator 
\begin{align*}
	\mathcal L=y^{\alpha_1}\mbox{Tr }\left(QD^2_x\right)+2y^{\frac{\alpha_1+\alpha_2}{2}}q\cdot \nabla_xD_y+\gamma y^{\alpha_2} D_{yy}+c y^{\alpha_2-1}D_y
\end{align*}
under Neumann boundary conditions. The methods use in these papers rely on tools from vector-valued harmonic analysis, Fourier multipliers and structure theory of Banach spaces.  We refer the reader also to  \cite{dong2020neumann, dong2020RMI, dong2020parabolic, dong2021weighted} and to  \cite{Vita2024higher,Vita2024schauder, dong-Vita2024schauder,Vita2023HarnackDini,Vita2021LiouvilleEven, Vita2021LiouvilleOdd,Vita2024Harnack} for related results with different methods, but without the powers $y^{\alpha_1}$, $y^{\alpha_2}$ ($\alpha_1=\alpha_2=0$) and with variable coefficients.

\medskip
This paper is devoted to complete the picture in this direction, by adding  the $x$-derivative  $y^{\frac{\alpha_1+\alpha_2}{2}-1}\left(d,\nabla_x\right)$ and the potential term  $by^{\alpha_2-2}$  and by studying $\mathcal L$ under Dirichlet or  oblique derivative boundary condition at $y=0$.
 Here we consider only constant matrices $Q$ and constant $q, \gamma$. The general case where $Q,q, \gamma$ are bounded and uniformly continuous is however straightforward and allows to treat  operators in smooth domains, whose degeneracy in the top order coefficients behaves like a power of the distance from the boundary. We shall treat these topics in a forthcoming paper.
 
  We also point out that our results seems to be new in the case of   oblique derivative boundary conditions (see Theorems \ref{Teo alpha Oblique} and \ref{complete-oblique}) when $\alpha_1\neq 0$ or $\alpha_2\neq 0$ (i.e. when the   powers $y^{\alpha_1}$, $y^{\alpha_2}$  appear in the operator)    while for Dirichlet boundary conditions (Theorems \ref{complete} and \ref{complete dirichlet alpha diff}) we improve the results in \cite{dong2024DivergenceAlpha} which are valid in the special case $\alpha_1=\alpha_2\in (0,2)$, $v=0$, $b=0$ but with variable coefficients.

\bigskip

\medskip

\medskip

The paper is organized as  follows. In Section \ref{Section main results} we collect the  hypotheses  we assume on $\mathcal L$ and present the main  results of the paper: here, in order to improve readability,  we restrict ourselves to the case $\alpha_1=\alpha_2$ as in  \eqref{defL alpha eq}, referring to Section  \ref{Consequenses} for the general case \eqref{defL}.
  
In Section \ref{Section Degenerate}, we exploit some elementary changes of variables, in a functional analytic setting, to reduce our operators to  simpler cases. 

In Section \ref{Weighted}  we collect the results we need concerning anisotropic weighted Sobolev spaces:  the main novelty is Section \ref{section oblique sobolev} where we introduce and study the Sobolev space $W^{2,p}_{v}(\alpha_1,\alpha_2,m)$ having oblique derivative boundary condition.
 This careful study is essential for characterizing the domain of  $\mathcal L$ when the $x$-derivative  $y^{\frac{\alpha_1+\alpha_2}{2}-1}\left(d,\nabla_x\right)$ is present and when $\alpha_1\neq \alpha_2$.

In Sections \ref{Sect DOm obliqu}  and \ref{Section DIrichlet alpha}, which are the core of the paper,  we prove generation results, maximal regularity and domain characterization for the operator $\mathcal L$, under respectively oblique derivative and Dirichlet boundary conditions, both  in  the case $\alpha_1=\alpha_2$.  

Finally, in Section \ref{Consequenses}, we extend the results to  general $\alpha_1$, $\alpha_2$.

\bigskip
\noindent\textbf{Notation.} For $N \ge 0$, $\R^{N+1}_+=\{(x,y): x \in \R^N, y>0\}$. For $m \in \R$ we consider the measure $y^m dx dy $ in $\R^{N+1}_+$ and  we write $L^p_m(\R_+^{N+1})$ for  $L^p(\R_+^{N+1}; y^m dx dy)$ and often only $L^p_m$ when $\R^{N+1}_+$ is understood. 
$\C^+=\{ \lambda \in \C: \Rp \lambda >0 \}$ and, for $|\theta| \leq \pi$, we denote by  $\Sigma_{\theta}$  the open sector $\{\lambda \in \C: \lambda \neq 0, \ |Arg (\lambda)| <\theta\}$.  
We denote by  $\alpha^+$ and $\alpha^-$ the positive and negative  part of a real  number, that is $\alpha^+=\max\{\alpha,0\}$,  $\alpha^-=-\min\{\alpha,0\}$.

We write often $(x,y)$ or $x\cdot y$ to denote the inner product of $\R^N$ and, for $A,B\in\R^{N,N}$ symmetric, $\mbox{Tr }\left(AB\right)=\sum_{i,j}a_{ij}b_{i,j}$. Moreover, if $\omega\in\R^N$, we write also  $\omega\otimes\omega\in\R^{N,N}$ to denote the matrix $\left(\omega_i\omega_j\right)_{i;j=1,\dots N}$; with this notation one has $\mbox{Tr }\left(\omega\otimes\omega  A\right)=\left(A\omega,\omega\right)$.

We use $B$ for the one-dimensional Bessel operator $D_{yy} +\frac cy D_y$ and $L$ for $D_{yy}+\frac{c}{y}D_y-\frac{b}{y^2}$. Here $c, b \in \R$ and both operators are defined on the half-line $(0, \infty)$.
\bigskip

\section{The main results and assumptions}\label{Section main results}
We consider first, for $b,c\in\R$, the 1d operators
\begin{align}\label{def 1-d op}
	L=D_{yy}+\frac{c}{y}D_y-\frac{b}{y^2}, \qquad B=D_{yy}+\frac{c}{y}D_y
\end{align}
on the  half line $\R_+=]0, \infty[$. Note that $B$ (which stands for Bessel) is nothing but $L$ when $b=0$. Often we write $L_y, B_y$ to indicate that they act with respect to the $y$ variable.



The equation $Lu=0$ has solutions $y^{-s_1}$, $y^{-s_2}$ where $s_1,s_2$ are the roots of the indicial equation $f(s)=-s^2+(c-1)s+b=0$ 

\begin{equation} \label{defs}
	s_1:=\frac{c-1}{2}-\sqrt{D},
	\quad
	s_2:=\frac{c-1}{2}+\sqrt{D}
\end{equation}
where
\bigskip
\begin{equation} \label{defD}
	D:=
	b+\left(\frac{c-1}{2}\right)^2.
\end{equation}

The above numbers are real if and only if $D \ge 0$. When $D<0$ the equation $u-Lu=f$ cannot have positive distributional solutions for certain positive $f$, see \cite{met-soba-spi3}. 
When $b=0$, then $\sqrt D=|c-1|/2$ and $s_1=0, s_2=c-1$ for $c \ge 1$ and $s_1=c-1, s_2=0$ for $c<1$.

\medskip

%

\medskip

We now introduce a $(N+1)$-d generalization of the operators \eqref{def 1-d op}. For reader's convenience we list below   the main hypothesis and  notation which we assume throughout  the whole manuscript.
 
\begin{assumption}\label{assumption}
	Let $v=(d,c)\in\R^{N+1}$ and let $A=\left(a_{ij}\right)\in \R^{N+1,N+1}$ be a symmetric and positive definite $(N+1)\times (N+1)$ matrix; we write $A$ as
\begin{align*}
A:=	\left(
\begin{array}{c|c}
	Q  & { q}^t \\[1ex] \hline
	q& \gamma
\end{array}\right)
\end{align*}
 where $Q\in \R^{N\times N}$,  $q=(q_1, \dots, q_N)\in\R^N$ and $\gamma=a_{N+1, N+1}>0$. Let $\alpha_1,\alpha_2\in\R$ such that 
 \begin{align*}
 	\alpha_2<2,\qquad   \alpha_2-\alpha_1<2.
 \end{align*} 
For $m \in \R$ we consider the measure $y^m dx dy $ in $\R^{N+1}_+$ and  we write $L^p_m(\R_+^{N+1})$ for  $L^p(\R_+^{N+1}; y^m dx dy)$ and often only $L^p_m$ when $\R^{N+1}_+$ is understood.
\end{assumption}

 \begin{defi}\label{Definition operator}
 	We consider the $(N+1)$-d  degenerate operator
\begin{align*}
	\mathcal L=y^{\alpha_1}\mbox{Tr }\left(QD^2_x\right)+2y^{\frac{\alpha_1+\alpha_2}{2}}q\cdot \nabla_xD_y+\gamma y^{\alpha_2} D_{yy}+y^{\frac{\alpha_1+\alpha_2}{2}-1}\left(d,\nabla_x\right)+cy^{\alpha_2-1}D_y-by^{\alpha_2-2}
\end{align*}
in the space $L^p_m=L^p_m(\R^{N+1}_+)$. Note that $\mathcal L$ can be written equivalently as
\begin{align}\label{def L transf}
	\mathcal L&=y^{\alpha_1}\mbox{Tr }\left(QD^2_x\right)+2y^{\frac{\alpha_1+\alpha_2}2}\left(q,\nabla_xD_y\right)+\gamma y^{\alpha_2}  L_y+y^{\frac{\alpha_1+\alpha_2}{2}-1}\left(d,\nabla_x\right)
\end{align}
where $L_y$ is the operator defined in \eqref{def 1-d op} with parameters $\frac b\gamma,\ \frac c\gamma$.
\medskip 

Note that in the special case $\alpha_1=\alpha_2=\alpha$, so the unique requirement is $\alpha<2$, the operator $\mathcal L$ takes the form
\begin{align}\label{def L transf alpha}
	\mathcal L&=y^{\alpha}\mbox{Tr }\left(AD^2\right)+y^{\alpha-1}\left(v,\nabla\right)-by^{\alpha-2}.
\end{align}

 \end{defi}
\begin{assumption}\label{Assumption D}
We always keep the assumption  $D\geq 0$ satisfied by the coefficients of the operator $ L_y$ in \eqref{def L transf}.
\end{assumption}


We study unique solvability of the problems 
\begin{align*}
	\lambda u-\mathcal L u=f,\qquad D_t v -\mathcal L v=g
\end{align*} 
in the spaces $L^p_m(\R^{N+1}_+)$ spaces under 
 Dirichlet or oblique derivative boundary conditions at $y=0$, and initial conditions in the parabolic case, together with the regularity of $u,v$. In the language of semigroup theory, we prove that $\mathcal L$ generates an analytic semigroup on $L^p_m$, characterize its domain and show that it has maximal regularity, which means that both $D_t v$ and $\mathcal L v$ have the same regularity as $g$. To improve  readability we recall the following definition.


%

\begin{defi}\label{defi Maximal regularity}
An analytic semigroup $(e^{t\mathcal A})_{t \ge0}$ on a Banach space $X$ with generator $\mathcal A$ has
{\it maximal regularity of type $L^q$} ($1<q<\infty$)
if for each $f\in L^q([0,T];X)$,   the following parabolic problem associated with $\mathcal A$
\begin{align}\label{Parabolic problem}
\left\{
\begin{array}{ll}    
	D_tu(t)-\mathcal Au(t)=f(t), \quad t>0, \\[1ex]
	u(0)=0 \end{array}\right.
\end{align}
has a unique solution $u\in W^{1,q}([0,T];X)\cap L^q([0,T];D(\mathcal A))$. This means that the mild solution of \eqref{Parabolic problem}, given by the variation of parameters formula
$$t\mapsto u(t)=\int_0^te^{(t-s)\mathcal A}f(s)\,ds,$$
  is indeed  a strong solution and has the best regularity one can expect.
\end{defi}
It is known that the property above does not depend on $1<q<\infty$ and $T>0$.
A characterization of maximal regularity is available in UMD Banach spaces, through the $\mathcal  R$-boundedness of the resolvent in a suitable sector $\omega+\Sigma_{\phi}$, with $\omega \in \R$ and $\phi>\pi/2$: this approach is widely described in \cite{KunstWeis} and in the new books \cite{WeisBook1}, \cite{WeisBook2}.\\

\medskip

Our main results are the following. We  refer to  Section \ref{Weighted}  for the definition of the weighted Sobolev spaces involved 
\begin{align*}
	W^{2,p}_{\mathcal N}\left(\alpha_1,\alpha_2,m\right),\qquad W^{2,p}_{v}\left(\alpha_1,\alpha_2,m\right),
\end{align*}
having respectively Neumann and oblique boundary conditions (see also Definition \ref{convention Wv Wn alpha}).\\ 

  For simplicity, here, we state the results only in the special  case $\alpha_1=\alpha_2=\alpha$, where 
\begin{align*}
	\mathcal L&=y^{\alpha}\mbox{Tr }\left(AD^2\right)+y^{\alpha-1}\left(v,\nabla\right)-by^{\alpha-2},
\end{align*}
 referring to Section \ref{Consequenses} for the general case having possibly different weights $y^{\alpha_1}$, $y^{\alpha_2}$, in front of the $x$ and $y$ derivatives.\\

We start by considering $b=0$ and    endow $\mathcal L$ with  the Neumann or  oblique derivative boundary conditions 
$$\lim_{y\to 0} D_y u=0\quad \text{(if $v=0$)},\qquad \lim_{y\to 0}y^{\frac c\gamma}\, v \cdot \nabla u=0\quad \text{(if $c\neq 0$)}.$$

We define accordingly (see Propositions \ref{neumann}, \ref{Trace D_yu in W}, \ref{Equiv Wv alpha eq})
\begin{align*}
	W^{2,p}(\alpha,\alpha,m)=&\left\{u\in W^{2,p}_{loc}(\R^{N+1}_+):\   y^{\alpha} D^2u,\ y^\frac{\alpha}{2} \nabla u,\ u\in L^p_m\right\}
\end{align*}
and
\begin{align*}
	W^{2,p}_{v}(\alpha, \alpha,m)&{=}\{u \in W^{2,p}(\alpha, \alpha,m):\ y^{\alpha-1}v\cdot \nabla u \in L^p_m\},\qquad (c\neq 0);\\[2ex]
	W^{2,p}_{v}(\alpha, \alpha,m){:=}W^{2,p}_{\mathcal N}(\alpha, \alpha,m)&{=}\{u \in W^{2,p}(\alpha, \alpha,m):\ y^{\alpha-1}D_yu \in L^p_m\},\qquad \hspace{1.8ex}(v= 0),
\end{align*}
\medskip 
(note that $W^{2,p}_{(0,c)}(\alpha, \alpha,m){=}W^{2,p}_{\mathcal N}(\alpha, \alpha,m)$).
\begin{teoNon}{\em (Theorems \ref{Teo alpha Oblique} and \ref{complete-oblique})} 
	Let $v=(d,c)\in\R^{N+1}$ with $d=0$ if $c=0$, and let $\alpha\in\R$ such that $\alpha<2$ and
	$$\alpha^{-} <\frac{m+1}p<\frac{c}{\gamma}+1-\alpha.$$   
	Then the operator
	\begin{align*}
		\mathcal L&=y^{\alpha}\mbox{Tr }\left(AD^2\right)+y^{\alpha-1}\left(v,\nabla\right)
	\end{align*}
	generates a bounded analytic semigroup  in $L^p_m$ which has maximal regularity. Moreover 
	\begin{align*}
		D(\mathcal L)=W^{2,p}_{v}\left(\alpha,\alpha,m\right)
	\end{align*} 
	and the set  $ \mathcal C_v$ defined in \eqref{defC oblique} is a core for $\mathcal L$.
\end{teoNon}

We then add the potential term $-by^{\alpha-2}$ and   endow $\mathcal L$ with the  Dirichlet boundary condition  (see  Corollary \ref{cor1})
\begin{align*}
	\lim_{y\to 0}y^{s_2} u=0 \quad (\text{if $D>0$}),\qquad 	\lim_{y\to 0}y^{s_2} u\in\C \quad (\text{if $D=0$}),
\end{align*}
where $D$ and  $s_{1,2}$ are defined in \eqref{defs}, \eqref{defD}.

\begin{teoNon}{\em (Theorems \ref{complete} and \ref{complete dirichlet alpha diff})}	Let $\alpha\in\R$ such that $\alpha<2$ and
	$$s_1+ \alpha^-<\frac{m+1}p<s_2+2-\alpha.$$   
	Then, under Assumptions   \ref{assumption} and \ref{Assumption D}, the operator
	\begin{align*}
		\mathcal L&=y^{\alpha}\mbox{Tr }\left(AD^2\right)+y^{\alpha-1}\left(v,\nabla\right)-by^{\alpha-2}
	\end{align*} 
	generates a bounded analytic semigroup  in $L^p_m$ which has maximal regularity. Moreover,
	\begin{equation*}
		D(\mathcal L)
		=y^{-s_1}W^{2,p}_{w}\left(\alpha,\alpha,m-s_1p\right),\qquad w=v-2s_1\left(q,\gamma\right).
	\end{equation*}	
	
	
\end{teoNon}

\medskip

The  maximal regularity of $\mathcal L$ stated in the Theorems above,  implies immediately  the following   result which we state, for simplicity, only in the case of oblique boundary conditions. The proof follows directly from the above theorems, Definition \ref{defi Maximal regularity} and standard semigroup theory.

\begin{cor}
		Let $v=(d,c)\in\R^{N+1}$ with $d=0$ if $c=0$, and let $\alpha\in\R$ such that $\alpha<2$ and
	$\alpha^{-} <\frac{m+1}p<\frac{c}{\gamma}+1-\alpha$.   
Let  us consider the operator 
$$\mathcal L=y^{\alpha}\mbox{Tr }\left(AD^2\right)+y^{\alpha-1}\left(v,\nabla\right)$$ 
	endowed with domain  $W^{2,p}_{v}\left(\alpha,\alpha,m\right):=W^{2,p}_{m,v}$. Then for each $1<q<\infty$, $T>0$ and $u_0 \in W^{2,p}_{m, v}$, $f\in L^q([0,T];L^p_m)$ the problem
	$$\begin{cases}
		\frac{\partial}{\partial t} u(t,x,y)-\mathcal Lu(t,x,y)=f(t,x,y), \quad t>0,\\[1ex]
		u(0,x,y)=u_0(x,y)
	\end{cases}$$
	admits a unique solution  $u\in W^{1,q}([0,T];L^p_m)\cap L^q([0,T];W^{2,p}_{m,v})$.
	
\end{cor}

\section{Degenerate operators and similarity transformations }\label{Section Degenerate}

In this section we consider the operator $\mathcal L$ defined in Definition \ref{Definition operator} namely

\begin{align*}	
	\mathcal L&=y^{\alpha_1}\mbox{Tr }\left(QD^2_x\right)+2y^{\frac{\alpha_1+\alpha_2}{2}}q\cdot \nabla_xD_y+\gamma y^{\alpha_2} D_{yy}+y^{\frac{\alpha_1+\alpha_2}{2}-1}\left(d,\nabla_x\right)+cy^{\alpha_2-1}D_y-by^{\alpha_2-2}
\end{align*}
which we often shorten by writing
\begin{align}\label{shorten L}	
	\mathcal L=y^{\alpha_1}\mbox{Tr }\left(QD^2_x\right)+2y^{\frac{\alpha_1+\alpha_2}2}\left(q,\nabla_xD_y\right)+\gamma y^{\alpha_2} L_y+y^{\frac{\alpha_1+\alpha_2}{2}-1}\left(d,\nabla_x\right)
\end{align}
where $L_y=D_yy+\frac{c/\gamma}{y}D_y-\frac{b/\gamma}{y^2}$ is the operator defined in \eqref{def 1-d op} with parameters $\frac b\gamma,\ \frac c\gamma$. \\

We investigate how this operator can be transformed  by means of  change of variables and multiplications.

For  $k,\beta \in\R$, $\beta\neq -1$ let
\begin{align}\label{Gen Kelvin def}
	T_{k,\beta\,}u(x,y)&:=|\beta+1|^{\frac 1 p}y^ku(x,y^{\beta+1}),\quad (x,y)\in\R^{N+1}_+.
\end{align}
Observe that
$$ T_{k,\beta\,}^{-1}=T_{-\frac{k}{\beta+1},-\frac{\beta}{\beta+1}\,}.$$

\begin{prop}\label{Isometry action der} Let $1\leq p\leq \infty$, $k,\beta \in\R$, $\beta\neq -1$. The following properties hold.
	\begin{itemize}
		\item[(i)] For every $m\in\R$,  $T_{k,\beta\,}$ maps isometrically  $L^p_{\tilde m}$ onto $L^p_m$  where 
		$$ \tilde m=\frac{m+kp-\beta}{\beta+1}.$$
		\item[(ii)] For every  $u\in W^{2,1}_{loc}\left(\R^{N+1}_+\right)$ one has
		\begin{itemize}
			\item[1.] $y^\alpha T_{k,\beta\,}u=T_{k,\beta\,}(y^{\frac{\alpha}{\beta+1}}u)$, for any $\alpha\in\R$;\medskip
			\item [2.] $D_{x_ix_j}(T_{k,\beta\,}u)=T_{k,\beta} \left(D_{x_ix_j} u\right)$, \quad $D_{x_i}(T_{k,\beta\,}u)=T_{k,\beta}\left(D_{x_i} u\right)$;\medskip
			\item[3.]  $D_y T_{k,\beta\,}u=T_{k,\beta\,}\left(ky^{-\frac 1 {\beta+1}}u+(\beta+1)y^{\frac{\beta}{\beta+1}}D_yu\right)$,
			\\[1ex] $D_{yy} (T_{k,\beta\,} u)=T_{k,\beta\,}\Big((\beta+1)^2y^{\frac{2\beta}{\beta+1}}D_{yy}u+(\beta+1)(2k+\beta)y^{\frac{\beta-1}{\beta+1}}D_y u+k(k-1)y^{-\frac{2}{\beta+1}}u\Big)$.\medskip
			\item[4.] $D_{xy} T_{k,\beta\,}u=T_{k,\beta\,}\left(ky^{-\frac 1 {\beta+1}}D_xu+(\beta+1)y^{\frac{\beta}{\beta+1}}D_{xy}u\right)$
		\end{itemize}
	\end{itemize}
\end{prop}{\sc{Proof.}} The proof of (i) follows after observing the Jacobian of $(x,y)\mapsto (x,y^{\beta+1})$ is $|1+\beta|y^{\beta}$. To prove (ii) one can easily observe that any $x$-derivatives commutes with $T_{k,\beta}$. Then we compute  
\begin{align*}
	D_y T_{k,\beta\,}u(x,y)=&|\beta+1|^{\frac 1 p}y^{k}\left(k\frac {u(x,y^{\beta+1})} y+(\beta+1)y^\beta D_y u(x,y^{\beta+1})\right)\\[1ex]
	=&T_{k,\beta\,}\left(ky^{-\frac 1 {\beta+1}}u+(\beta+1)y^{\frac{\beta}{\beta+1}}D_yu\right)
\end{align*}
and similarly
\begin{align*}
	D_{yy} T_{k,\beta\,} u(x,y)=&T_{k,\beta\,}\Big((\beta+1)^2y^{\frac{2\beta}{\beta+1}}D_{yy}u+(\beta+1)(2k+\beta)y^{\frac{\beta-1}{\beta+1}}D_y u+k(k-1)y^{-\frac{2}{\beta+1}}u\Big).
\end{align*}
\qed

\begin{prop}\label{Isometry action}  Let 
	$T_{k,\beta\,}$ be the isometry above defined and let
	\begin{align*}	
		\mathcal L&=y^{\alpha_1}\mbox{Tr }\left(QD^2_x\right)+2y^{\frac{\alpha_1+\alpha_2}{2}}q\cdot \nabla_xD_y+\gamma y^{\alpha_2}D_{yy}+y^{\frac{\alpha_1+\alpha_2}{2}-1}\left(d,\nabla_x\right)+c y^{\alpha_2-1}D_y-by^{\alpha_2-2}.
	\end{align*}
   The following properties hold.
	
	\begin{itemize}
		\item[(i)] For every  $u\in W^{2,1}_{loc}\left(\R^{N+1}_+\right)$ one has
		\begin{align*}
			&T_{k,\beta\,}^{-1} \Big(\mathcal L		\Big)T_{k,\beta\,}u=\tilde{\mathcal L}\,u
		\end{align*}
	where 
	\begin{align*}
	\tilde{\mathcal L}=y^{\tilde{\alpha}_1}\mbox{Tr }\left(QD^2_x\right)+2y^{\frac{{\tilde \alpha}_1+{\tilde \alpha}_2}{2}}\tilde q\cdot \nabla_xD_y+y^{{\tilde \alpha}_2}\tilde \gamma D_{yy}+y^{\frac{{\tilde \alpha}_1+{\tilde \alpha}_2}{2}-1}\left(\tilde d,\nabla_x\right)+y^{{\tilde \alpha}_2-1}\tilde cD_y-\tilde by^{{\tilde \alpha}_2-2}
	\end{align*}
 is the operator defined as in Definition \ref{Definition operator}	but with parameters given by $$\tilde\alpha_1=\frac{\alpha_1}{\beta+1},\qquad\tilde\alpha_2=\frac{\alpha_2+2\beta}{\beta+1}$$
  and  
		\begin{align}\label{tilde b}
		\nonumber	\tilde q=(\beta+1)q,\qquad\quad  \tilde \gamma =(\beta+1)^2\gamma,\qquad \quad 
\tilde d=2kq+d,\\[1ex]
\tilde c=(\beta+1)\left(c+(2k+\beta)\gamma\right),\qquad \quad 	\tilde b=b-k\left(c+(k-1)\gamma\right).
		\end{align}
	\item[(ii)]	In particular choosing $\beta=\frac{\alpha_1-\alpha_2}2$ and setting $\tilde \alpha=\frac{2\alpha_1}{\alpha_1-\alpha_2+2}$ one has $\tilde \alpha_1=\tilde \alpha_2=\tilde\alpha$ and, for every  $u\in W^{2,1}_{loc}\left(\R^{N+1}_+\right)$,
		\begin{align*}
			T_{k,\beta\,}^{-1} \Big(\mathcal L\Big)T_{k,\beta\,}u
			&=y^{\tilde \alpha}\mbox{Tr }\left(\tilde AD^2u\right)+y^{\tilde \alpha-1}\left(\tilde v,\nabla u \right)-\tilde by^{\alpha-2}u.
		\end{align*}
	where $\tilde A=	\left(
	\begin{array}{c|c}
		Q  & { \tilde q}^t \\[1ex] \hline
		\tilde q& \tilde \gamma
	\end{array}\right)$ and $\tilde v=(\tilde d, \tilde c)$.
	\end{itemize}	
\end{prop} 
{\sc{Proof.}} The proof follows  after a tedious straightforward computation using Proposition \ref{Isometry action der}.
\qed

If in the above Proposition we write both operator $\mathcal L$, $\tilde{\mathcal L}$ in the compact form \eqref{shorten L}, we have the following result.
\begin{cor}\label{Isometry action2}  Let 
	$T_{k,\beta\,}$ be the isometry above defined and let 
	\begin{align*}	
		\mathcal L=y^{\alpha_1}\mbox{Tr }\left(QD^2_x\right)+2y^{\frac{\alpha_1+\alpha_2}2}\left(q,\nabla_xD_y\right)+\gamma y^{\alpha_2} L_y+y^{\frac{\alpha_1+\alpha_2}{2}-1}\left(d,\nabla_x\right)
	\end{align*}
	where $L_y=D_{yy}+\frac{c/\gamma}{y}D_y-\frac{b/\gamma}{y^2}$. 
	 Then
	
	\begin{itemize}
		\item[(i)] For every  $u\in W^{2,1}_{loc}\left(\R^{N+1}_+\right)$ one has
		\begin{align*}
			T_{k,\beta\,}^{-1} \Big(\mathcal L\Big)T_{k,\beta\,}u&=y^{\tilde\alpha_1}\mbox{Tr }\left(QD^2_xu\right)+2y^{\frac{\tilde\alpha_1+\tilde\alpha_2}2}\left(\tilde q,\nabla_xD_yu\right)+y^{\tilde\alpha_2}\tilde \gamma  \tilde L_yu+y^{\frac{\tilde\alpha_1+\tilde\alpha_2}{2}-1}\left(\tilde d,\nabla_x u\right)
		\end{align*}
		where  $\tilde { L}_y=D_{yy}+\frac{\tilde c/\tilde\gamma }{y}D_y-\frac{\tilde b/\tilde\gamma}{y^2}$.
%
		\item[(ii)] The discriminant $\tilde D$ and the parameters $\tilde s_{1,2}$ of $\tilde L_y$ defined as in \eqref{defD}, \eqref{defs}  are related to  those of $L_y$ by
		\begin{align}\label{tilde D gamma}
			\tilde D&=\frac{D}{(\beta+1)^2},
		\end{align}
		and
		\begin{align}\label{tilde s gamma}
			\tilde s_{1,2}=\frac{s_{1,2}+k}{\beta+1} \quad  (\beta+1>0), \qquad
			\tilde s_{1,2}&=\frac{s_{2,1}+k}{\beta+1}\quad  (\beta+1<0).
		\end{align}
	\end{itemize}	
\end{cor} 
{\sc{Proof.}} The first claim is simply a reformulation of (i) of Proposition \ref{Isometry action}. (ii) then follows directly by \eqref{tilde b} and Definitions \eqref{defs}, \eqref{defD}.

\qed

We define now for $\omega\in\R^N$ and $\beta\neq -1$ the following isometry of $L^p_m$
\begin{align}\label{Tran def}
	S_{\beta,\omega}\, u(x,y)&:=u\left(x+\omega y^{\beta+1},y\right),\quad (x,y)\in\R_+^{N+1}.
\end{align}

\begin{prop}\label{Isometry shift} Let  $\omega\in\R^N$ and $\beta\neq -1$. Then for every $m\in\R$,  $ S_{\beta,\omega}$ is an isometry of   $L^p_{m}$ and for every  $u\in W^{2,1}_{loc}\left(\R^{N+1}_+\right)$ one has
	\begin{itemize}
		\item[1.] $y^\alpha S_{\beta,\omega} u={S_{\beta,\omega}}(y^{\alpha}u)$, for any $\alpha\in\R$;\medskip
		\item [2.] $D_{x_ix_j}({S_{\beta,\omega}}u)={S_{\beta,\omega}} \left(D_{x_ix_j} u\right)$, \quad $D_{x_i}({S_{\beta,\omega}}u)={S_{\beta,\omega}}\left(D_{x_i} u\right)$;\medskip
		\item[3.]  $D_y {S_{\beta,\omega}}u={S_{\beta,\omega}}\Big((\beta+1)y^\beta (\nabla_x u,\omega)+D_yu\Big)$,
		\\[1ex] $D_{yy} ({S_{\beta,\omega}} u)={S_{\beta,\omega}}\Big((\beta+1)^2y^{2\beta} (D^2_x u\cdot\omega,\omega)+2(\beta+1)y^\beta (\nabla_x D_y u, \omega)$
		\begin{flushright}
			\vspace{-2ex}$
			+\beta (\beta+1)y^{\beta-1}(\nabla_x u,\omega)+ D_{yy}u\Big)$;
		\end{flushright}
		\item[4.] $\nabla_xD_{y} {S_{\beta,\omega}}u={S_{\beta,\omega}}\Big((\beta+1)y^\beta D^2_x u\cdot \omega+\nabla_xD_yu\Big)$.
	\end{itemize}
\end{prop}
{\sc{Proof.}} The proof follows after a straightforward computation.\qed

\begin{prop}\label{Isometry action 2}  Let 
	${S_{\beta,\omega}}$ be the isometry above defined and let $\mathcal L$ the operator defined in  \eqref{def L transf} with $L_y=D_{yy}+\frac{c/\gamma}{y}D_y-\frac{b/\gamma}{y^2}$.   Then for every  $u\in W^{2,1}_{loc}\left(\R^{N+1}_+\right)$ one has
	\begin{align*}
		&{S_{\beta,\omega}}^{-1} \Big(y^{\alpha_1}\mbox{Tr }\left(QD^2_x\right)+2y^{\frac{\alpha_1+\alpha_2}2}\left(q,\nabla_xD_y\right)+\gamma y^{\alpha_2} L_y+y^{\frac{\alpha_1+\alpha_2}{2}-1}\left(d,\nabla_x\right)\Big){S_{\beta,\omega}}u\\[1ex]
		&=y^{\alpha_1}\mbox{Tr }\left(QD^2_xu\right)+2(\beta+1)y^{\frac{\alpha_1+\alpha_2+2\beta}2}\mbox{Tr }(q\otimes\omega D^2_xu)+\gamma (\beta+1)^2y^{\alpha_2+2\beta}\mbox{Tr }\left(\omega\otimes\omega D^2_xu\right)\\[1ex]
		&+2y^{\frac{\alpha_1+\alpha_2}2}\left(q,\nabla_xD_yu\right)+2\gamma (\beta+1)y^{\alpha_2+\beta}(\omega,\nabla_xD_y u)\\[1ex]
		&+\gamma y^{\alpha_2}D_{yy}u+cy^{\alpha_2-1}D_yu-by^{\alpha_2-2}u\\[1ex]
		&+(c+\gamma\beta)(\beta+1)y^{\alpha_2+\beta-1}(\omega, \nabla_xu)+y^{\frac{\alpha_1+\alpha_2}2-1}(d,\nabla_xu).
	\end{align*}
\end{prop} 
{\sc{Proof.}} Using Proposition \ref{Isometry action der}, one has 
\begin{itemize}
	\item ${S_{\beta,\omega}}^{-1} \Big(y^{\alpha_1}\mbox{Tr }\left(QD^2_xu\right)\Big){S_{\beta,\omega}}u=y^{\alpha_1}\mbox{Tr }\left(QD^2_xu\right)$;
	\item ${S_{\beta,\omega}}^{-1} \Big(2y^{\frac{\alpha_1+\alpha_2}2}\left(q,\nabla_xD_y\right)\Big){S_{\beta,\omega}}u=2(\beta+1)y^{\frac{\alpha_1+\alpha_2+2\beta}2}\mbox{Tr }(q\otimes\omega D^2_xu)+2y^{\frac{\alpha_1+\alpha_2}2}\left(q,\nabla_xD_yu\right)$;
	\item ${S_{\beta,\omega}}^{-1} \Big(y^{\alpha_2}\gamma D_{yy}\Big){S_{\beta,\omega}}u=\gamma (\beta+1)^2y^{\alpha_2+2\beta}\mbox{Tr }\left(\omega\otimes\omega D^2_xu\right)+2\gamma (\beta+1)y^{\alpha_2+\beta}(\omega,\nabla_xD_y u)$
	\begin{flushright}
		\vspace{-2ex}	$+\gamma y^{\alpha_2}D_{yy}u+\gamma \beta(\beta+1)y^{\alpha_2+\beta-1}(\omega,\nabla_x u)$;
	\end{flushright}
	\item ${S_{\beta,\omega}}^{-1} \Big(cy^{\alpha_2-1}D_y\Big){S_{\beta,\omega}}u=c(\beta+1)y^{\alpha_2+\beta-1}(\omega, \nabla_xu)+cy^{\alpha_2-1}D_yu$;
	\item ${S_{\beta,\omega}}^{-1} \Big(y^{\frac{\alpha_1+\alpha_2}{2}-1}\left(d,\nabla_x\right)\Big){S_{\beta,\omega}}u=y^{\frac{\alpha_1+\alpha_2}2-1}(d,\nabla_xu)$.
\end{itemize}.
 The required claim then  follows after a straightforward computation.
\qed

When $\alpha_2-\alpha_1\neq 2$ we can specialize the previous relation by choosing $\beta=\frac{\alpha_1-\alpha_2}2$.

\begin{cor}\label{Isometry action 2 spe}  Let  $\omega\in\R^N$ and let ${S_{\beta,\omega}}$ be the isometry defined in \eqref{Tran def} with  $\beta=\frac{\alpha_1-\alpha_2}2\neq -1$. Let $\mathcal L$ the operator defined in  \eqref{def L transf} with $L_y=D_{yy}+\frac{c/\gamma}{y}D_y-\frac{b/\gamma}{y^2}$.   Then one has
	\begin{align*}
		{S_{\beta,\omega}}^{-1} &\Big(y^{\alpha_1}\mbox{Tr }\left(QD^2_x\right)+2y^{\frac{\alpha_1+\alpha_2}2}\left(q,\nabla_xD_y\right)+\gamma y^{\alpha_2} L_y+y^{\frac{\alpha_1+\alpha_2}{2}-1}\left(d,\nabla_x\right)\Big){S_{\beta,\omega}}u
		\\[1ex]
		&=y^{\alpha_1}\mbox{Tr }\left(\tilde QD^2_xu\right)+2y^{\frac{\alpha_1+\alpha_2}2}\left(\tilde q,\nabla_xD_yu\right)+\gamma y^{\alpha_2} L_yu+y^{\frac{\alpha_1+\alpha_2}{2}-1}\left(\tilde d,\nabla_xu\right)
	\end{align*}
	where 
	\begin{align*}
		\tilde Q=Q+2(\beta+1)&q\otimes \omega+\gamma(\beta+1)^2\omega\otimes\omega,\qquad \tilde q=q+\gamma(\beta+1)\omega,\\[1ex] &\tilde d=d+(c+\gamma\beta)(\beta+1)\omega
	\end{align*}	
\end{cor} 
{\sc{Proof.}} The proof follows by specializing Proposition \ref{Isometry action 2} to $\beta=\frac{\alpha_1-\alpha_2}2\neq -1$. 
\qed

\section{Anisotropic weighted Sobolev spaces} \label{Weighted} 
Let $p>1$, $m, \alpha_1,\alpha_2 \in \R$ such that
\begin{align}\label{param sobolev}
	\alpha_2<2,\qquad \alpha_2-\alpha_1<2,\qquad \alpha_1^{-} <\frac{m+1}p.
\end{align}
In order to describe the domain of  the operator \begin{align}\label{def operator sobolev section}
	\mathcal L&=y^{\alpha_1}\mbox{Tr }\left(QD^2_x\right)+2y^{\frac{\alpha_1+\alpha_2}2}\left(q,\nabla_xD_y\right)+\gamma y^{\alpha_2} L_y+y^{\frac{\alpha_1+\alpha_2}{2}-1}\left(d,\nabla_x\right)
\end{align}
we collect in this section  the  results we need  about suitable anisotropic weighted Sobolev spaces. 

The main novelty is Section \ref{section oblique sobolev} where we introduce the Sobolev space $W^{2,p}_{v}(\alpha_1,\alpha_2,m)$ having oblique derivative  boundary condition: here, in order to improve readability and although not essential,  we treat    separately
 the cases $\alpha_1=\alpha_2$ and $\alpha_1\neq\alpha_2$.
 
  Besides this we also briefly  recall the main properties about the spaces $W^{2,p}_{\mathcal N}(\alpha_1,\alpha_2,m)$ and  $W^{2,p}_{\mathcal R}(\alpha_1,\alpha_2,m)$ having respectively Neumann and Dirichlet boundary conditions,  referring to \cite{MNS-Sobolev} for  further details  and all the relative proofs (also outside the above range of parameters \eqref{param sobolev}).
  
   We also clarify, in Propositions \ref{Hardy Rellich Sob} and \ref{RN}, the relation between the three spaces $W^{2,p}_{\mathcal N}(\alpha_1,\alpha_2,m)$, $W^{2,p}_{v}(\alpha_1,\alpha_2,m)$ and  $W^{2,p}_{\mathcal R}(\alpha_1,\alpha_2,m)$.

\subsection{The space $W^{2,p}_{\mathcal N}(\alpha_1,\alpha_2,m)$}\label{Section neumann sobolev}
We start by defining the Sobolev space
\begin{align*}
	W^{2,p}(\alpha_1,\alpha_2,m)\stackrel{def}{=}&\left\{u\in W^{2,p}_{loc}(\R^{N+1}_+):\ u,\  y^{\alpha_1} D_{x_ix_j}u,\ y^\frac{\alpha_1}{2} D_{x_i}u,  \right.\\[1ex]
	&\left.\hspace{20ex}y^{\alpha_2}D_{yy}u,\ y^{\frac{\alpha_2}{2}}D_{y}u,\,y^\frac{\alpha_1+\alpha_2}{2} D_{y}\nabla_x u\in L^p_m\right\}
\end{align*}
which is a Banach space equipped with the norm
\begin{align*}
	\|u\|_{W^{2,p}(\alpha_1,\alpha_2,m)}\stackrel{def}{=}&\|u\|_{L^p_m}+\sum_{i,j=1}^n\|y^{\alpha_1} D_{x_ix_j}u\|_{L^p_m}+\sum_{i=1}^n\|y^{\frac{\alpha_1}2} D_{x_i}u\|_{L^p_m}\\[1ex]
	&+\|y^{\alpha_2}D_{yy}u\|_{L^p_m}+\|y^{\frac{\alpha_2}{2}}D_{y}u\|_{L^p_m}+\|y^\frac{\alpha_1+\alpha_2}{2} D_{y}\nabla_x u\|_{L^p_m}.
\end{align*}
Next we add different boundary conditions for $y=0$. 

We add a Neumann boundary condition for $y=0$  in the form $y^{\alpha_2-1}D_yu\in L^p_m$ and set
\begin{align*}
	W^{2,p}_{\mathcal N}(\alpha_1,\alpha_2,m)\stackrel{def}{=}\{u \in W^{2,p}(\alpha_1,\alpha_2,m):\  y^{\alpha_2-1}D_yu\ \in L^p_m\}
\end{align*}
with the norm
$$
\|u\|_{W^{2,p}_{\mathcal N}(\alpha_1,\alpha_2,m)}\stackrel{def}{=}\|u\|_{W^{2,p}(\alpha_1,\alpha_2,m)}+\|y^{\alpha_2-1}D_yu\|_{ L^p_m}.
$$
\begin{os}
	We remark  that, in the range of parameters \eqref{param sobolev}, the condition on the mixed derivatives in the definition of  $W^{2,p}_{\mathcal{N}}\left(\alpha_1,\alpha_2,m\right)$ can be discarded, without any loss of generality, since by 
	\cite[Proposition 4.1]{MNS-Sobolev} and \cite[Theorem 7.1]{MNS-CompleteDegenerate}
	one has  for every $u \in W^{2,p}_{\mathcal{N}}\left(\alpha_1,\alpha_2,m\right)$  
	\begin{align*}
		\|y^\frac{\alpha_1+\alpha_2}{2} D_{y}\nabla_x u \|_{ L^p_m} \leq& C \Big[\|u\|_{L^p_m}+\|y^{\alpha_1} D^2_{x}u\|_{L^p_m}+\|y^{\frac{\alpha_1}2} \nabla_{x}u\|_{L^p_m}\Big.
		\\[1ex]\Big.&\quad+\|y^{\alpha_2}D_{yy}u\|_{L^p_m}+\|y^{\frac{\alpha_2}{2}}D_{y}u\|_{L^p_m}+\|y^{\alpha_2-1}D_yu\|_{ L^p_m}\Big].
	\end{align*}
\end{os}

\begin{os}\label{Os Sob 1-d}
	With obvious changes we consider also the analogous Sobolev spaces $W^{2,p}(\alpha,m)$ and $W^{2,p}_{\cal N}(\alpha, m)$ on $\R_+$. 
	For example  we have
	$$W^{2,p}_{\mathcal N}(\alpha,m)=\left\{u\in W^{2,p}_{loc}(\R_+):\ u,\    y^{\alpha}D_{yy}u,\ y^{\frac{\alpha}{2}}D_{y}u,\ y^{\alpha-1}D_{y}u\in L^p_m\right\}.$$
	All the results of this section will be valid also in $\R_+$ changing (when it appears) the condition $\alpha_1^{-} <\frac{m+1}p$  to $0<\frac{m+1}p$.
\end{os}

The next result clarifies in which sense the condition $y^{\alpha_2-1}D_y u \in L^p_m$ is a Neumann boundary condition.

\begin{prop}{\em \cite[Proposition 4.3]{MNS-Sobolev}} \label{neumann} The following assertions hold.
	\begin{itemize} 
		\item[(i)] If $\frac{m+1}{p} >1-\alpha_2$, then $W^{2,p}_{\mathcal N}(\alpha_1, \alpha_2, m)=W^{2,p}(\alpha_1, \alpha_2, m)$.
		\item[(ii)] If $\frac{m+1}{p} <1-\alpha_2$, then $$W^{2,p}_{\mathcal N}(\alpha_1, \alpha_2, m)=\{u \in W^{2,p}(\alpha_1, \alpha_2, m): \lim_{y \to 0}D_yu(x,y)=0\ {\rm for\ a.e.\   x \in \R^N }\}.$$
	\end{itemize}
	In both cases (i) and (ii), the norm of $W^{2,p}_{\mathcal N}(\alpha_1, \alpha_2, m)$ is equivalent to that of $W^{2,p}(\alpha_1, \alpha_2, m)$.
\end{prop}

\medskip 

The next results  show the density  of smooth functions in $W^{2,p}_{\mathcal N}(\alpha_1,\alpha_2,m)$. Let
\begin{equation} \label{defC}
	\mathcal{C}:=\left \{u \in C_c^\infty \left(\R^N\times[0, \infty)\right), \ D_y u(x,y)=0\  {\rm for} \ y \leq \delta\ {\rm  and \ some\ } \delta>0\right \},
\end{equation}
its one dimensional version 
\begin{equation} \label {defDcore}
	\mathcal{D}=\left \{u \in C_c^\infty ([0, \infty)), \ D_y u(y)=0\  {\rm for} \ y \leq \delta\ {\rm  and \ some\ } \delta>0\right \}
\end{equation}
and finally (finite sums below)
$$C_c^\infty (\R^{N})\otimes\mathcal D=\left\{u(x,y)=\sum_i u_i(x)v_i(y), \  u_i \in C_c^\infty (\R^N), \  v_i \in \cal D \right \}\subset \mathcal C.$$
\begin{teo} {\em \cite[Theorem 4.9]{MNS-Sobolev}}\label{core gen}
	$C_c^\infty (\R^{N})\otimes\mathcal D$
 is dense in $W^{2,p}_{\mathcal N}(\alpha_1,\alpha_2,m)$.
\end{teo}

Note that the condition $(m+1)/p>\alpha_1^-$, or $m+1>0$ and $(m+1)/p+\alpha_1>0$, is necessary for the inclusion  $\mathcal C\subset W^{2,p}_{\mathcal N}(\alpha_1,\alpha_2,m)$.
\medskip


We provide  an equivalent description of $W^{2,p}_{\mathcal N}(\alpha_1, \alpha_2, m)$, adapted to the degenerate  operator $B_y=D_{yy}+cy^{-1}D_y$. In the first formulation we shows that the  Neumann boundary condition  in the integral  form $y^{\alpha_2-1}D_yu\in L^p_m$ is actually equivalent to the trace condition $\ds \lim_{y\to 0}y^c D_yu=0$. 
\begin{prop}{\em \cite[Proposition 4.6]{MNS-Sobolev}}\label{Trace D_yu in W}
	Let   $c\in\R$ and $\frac{m+1}{p}<c+1-\alpha_2$.  Then setting  $B_y=D_{yy}+cy^{-1}D_y$ one has 
	\begin{align*} 
		W^{2,p}_{\mathcal N}(\alpha_1, \alpha_2, m)=&\left\{u \in  W^{2,p}_{loc}(\R^{N+1}_+): u,\ y^{\alpha_1}\Delta_xu,\ y^{\alpha_2}B_y u \in L^p_m\text{\;\;and\;\;}\lim_{y\to 0}y^c D_yu=0\right\}
	\end{align*}
	and the norms $\|u\|_{W^{2,p}_{\mathcal N}(\alpha_1,\alpha_2,m)}$ and $$\|u\|_{L^p_m}+\|y^{\alpha_1}\Delta_x u\|_{L^p_m}+\|y^{\alpha_2}B_y u\|_{L^p_m}$$ are equivalent on $W^{2,p}_{\mathcal N}(\alpha_1, \alpha_2, m)$.
	Finally,  when $0<\frac{m+1}p\leq c-1$ then 
	\begin{align*} 
		W^{2,p}_{\mathcal N}(\alpha_1, \alpha_2, m)=&\left\{u \in  W^{2,p}_{loc}(\R^{N+1}_+): u,\ y^{\alpha_1}\Delta_xu,   y^{\alpha_2}B_y u \in L^p_m\right\}.
	\end{align*}
\end{prop}
The following equivalent description of $W^{2,p}_{\mathcal N}(\alpha_1, \alpha_2, m)$ involves a Dirichlet, rather than Neumann, boundary condition,  in a certain range of parameters.
\begin{prop}{\em \cite[Proposition 4.7]{MNS-Sobolev}}\label{trace u in W op}
	Let   $c\geq 1$ and $\frac{m+1}{p}<c+1-\alpha_2$.  The following properties hold.
	\begin{itemize}
		\item[(i)] If $c>1$ then 
		\begin{align*}
			W^{2,p}_{\mathcal N}(\alpha_1, \alpha_2, m)=&\left\{u \in  W^{2,p}_{loc}(\R^{N+1}_+): u,\ y^{\alpha_1}\Delta_xu,\ y^{\alpha_2}B_y \in L^p_m\text{\;and\;}\lim_{y\to 0}y^{c-1} u=0\right\}.
		\end{align*}
		\item[(ii)] If $c=1$ then 
		\begin{align*}
			W^{2,p}_{\mathcal N}(\alpha_1, \alpha_2, m)=&\left\{u \in  W^{2,p}_{loc}(\R^{N+1}_+): u,\ y^{\alpha_1}\Delta_xu, \ y^{\alpha_2}B_y u \in L^p_m\text{\;and\;}\lim_{y\to 0} u(x,y)\in \C\right\}.
		\end{align*}
	\end{itemize} 
\end{prop}
\subsection{The space $W^{2,p}_{v}(\alpha_1,\alpha_2,m)$}\label{section oblique sobolev}
\medskip 
Let  $v=(d,c)\in\R^{N+1}$, with $d\in\R^N$ and $c \neq 0$. We impose now  a weighted oblique derivative boundary condition  $\left(y^{\frac{\alpha_1-\alpha_2}2} d \cdot \nabla_x u+cD_yu\right)(x,0)=0$ in the integral form (see Proposition \ref{Equiv Wv alpha diff})
\begin{align}\label{eq bound cond alpha diff}
	y^{\frac{\alpha_1+\alpha_2}2-1} d \cdot \nabla_x u+cy^{\alpha_2-1}D_yu\in L^p_m.
\end{align}
For reader's convenience, although not necessary,   we treat    separately
the simpler case $\alpha_1=\alpha_2$ and the case  $\alpha_1\neq\alpha_2$, where some complications occur due to the different weights $y^{\alpha_1}$, $y^{\alpha_2}$ which appear in the $x$ and $y$ directions.

\subsubsection{The case $\alpha_1=\alpha_2$: the space $W^{2,p}_{v}(\alpha,\alpha,m)$}\label{section Wv alpha eq}
We start by the case $\alpha_1=\alpha_2:=\alpha$ where the condition above reads as 
$$y^{\alpha-1}v\cdot\nabla u=y^{\alpha-1}\left( d \cdot \nabla_x u+cD_yu\right)\in L^p_m.$$
We define accordingly
\begin{align}\label{Definition W_v alpha}
	W^{2,p}_{v}(\alpha, \alpha,m)\stackrel{def}{=}\{u \in W^{2,p}(\alpha, \alpha,m):\ y^{\alpha-1}v\cdot \nabla u \in L^p_m\}
\end{align}
with the norm
$$
\|u\|_{W^{2,p}_{v}(\alpha,\alpha,m)}\stackrel{def}{=}\|u\|_{W^{2,p}(\alpha,\alpha,m)}+\| y^{\alpha-1}v\cdot \nabla u\|_{ L^p_m}.
$$
In particular, when $d=0$, one has $W^{2,p}_{\mathcal N}(\alpha,\alpha,m)=W^{2,p}_{(0,c)}(\alpha,\alpha,m)$. This justify the following definition.
\begin{defi}
	To shorten some statements we also write
\begin{align*}
	W^{2,p}_{(0,0)}(\alpha,\alpha,m)\stackrel{def}{=}W^{2,p}_{\mathcal N}(\alpha,\alpha,m).
\end{align*}
With this notation $W^{2,p}_{v}(\alpha,\alpha,m)$ is well defined  for any $v=(d,c)\in\R^{N+1}$ such that $d=0$ if $c=0$. Note that  \begin{align*}
	W^{2,p}_{v}(\alpha,\alpha,m)=W^{2,p}_{\mathcal N}(\alpha,\alpha,m)\qquad \text{when\quad  $c=0$\quad  or\quad  $d=0$}.
\end{align*}
Under this identification, we also define
\begin{align} \label{defC oblique alpha}
	\mathcal{C}_v
	&:=\left \{u \in C_c^\infty \left(\R^N\times[0, \infty)\right), \ \left( v\cdot\nabla u\right)(x,y)=0\  {\rm for} \ y \leq \delta\ {\rm  and \ some\ } \delta>0\right \}
\end{align}
with the convention that $\mathcal{C}_{v}:=\mathcal C$ when $v=(0,0)$ and 	$\mathcal{C}$ is the set defined in \eqref{defC}. 
\end{defi}

\medskip 

In what follows we clarify the relation  between the spaces $W^{2,p}_{\mathcal N}(\alpha,\alpha,m)$ and   $W^{2,p}_{v}(\alpha,\alpha,m)$. The next Proposition shows that  $W^{2,p}_{v}(\alpha, \alpha,m)$ is related to $W^{2,p}_{\mathcal N}(\alpha,\alpha,m)$ by means of the isometry $S_{0,\omega}$ of $L^p_m$  defined in \eqref{Tran def} with  $\beta=0$ and  $\omega=-\frac d {c}$, namely
\begin{align*}
	S_{0,-\frac d {c}}\, u(x,y)&:=u\left(x-\frac d {c} y,y\right),\quad (x,y)\in\R_+^{N+1}.
\end{align*}

\begin{prop}\label{Equi obli Neumm sob alpha}
	One has 	
	\begin{align}\label{Equi obli Neumm sob alpha eq}
		S_{0,-\frac d c}\,\left(W^{2,p}_{\mathcal N}(\alpha, \alpha,m)\right)	=W^{2,p}_{v}(\alpha, \alpha,m)
	\end{align}
	In particular the set		$\mathcal{C}_v$ defined in \eqref{defC oblique alpha} 	is dense in $W^{2,p}_{v}(\alpha, \alpha,m)$.
\end{prop}
{\sc{Proof.}} Let $u\in W^{2,1}_{loc}\left(\R^{N+1}_+\right)$  and let us set $\tilde u={S_{0,-\frac d c}}u$. Then  by Proposition \ref{Isometry shift} one has 
\begin{itemize}
	\item [1.] $y^{\alpha}D_{x_ix_j}\tilde u={S_{0,-\frac d c}} \left(y^{\alpha}D_{x_ix_j} u\right)$, \quad $y^{\frac{\alpha}2}D_{x_i}\tilde u={S_{0,-\frac d c}}\left(y^{\frac{\alpha}2}D_{x_i} u\right)$,\\[1ex] $y^{\alpha-1}D_{x_i}\tilde u={S_{0,-\frac d c}}\left(y^{\alpha-1}D_{x_i} u\right)$;\medskip
	\item[2.]  $y^{\frac{\alpha}2}D_y \tilde u={S_{0,-\frac d c}}\Big(-\frac{1}cy^{\frac{\alpha}2} (\nabla_x u,d)+y^{\frac{\alpha}2}D_yu\Big)$,	\\[1ex]
	$y^{\alpha-1}D_y \tilde u={S_{0,-\frac d c}}\Big(-\frac 1 cy^{\alpha-1} (\nabla_x u,d)+y^{\alpha-1}D_yu\Big)$;\medskip
	\item[3.]	 $y^{\alpha_2}D_{yy} \tilde u={S_{0,-\frac d c}}\Big(\frac{1}{c^2}y^{\alpha} (D^2_x u\cdot d,d)-\frac 1 c2y^{\alpha} (\nabla_x D_y u, d)+ y^{\alpha}D_{yy}u\Big)$
	
	\item[4.] $y^{\alpha}\nabla_xD_{y} \tilde u={S_{0,-\frac d c}}\Big(-\frac 1 cy^{\alpha} D^2_x u\cdot d+y^{\alpha}\nabla_xD_yu\Big)$.
\end{itemize}
In particular 	$$y^{\alpha-1} d \cdot \nabla_x \tilde u+cy^{\alpha-1}D_y\tilde u ={S_{0,-\frac d c}}\Big(cy^{\alpha-1} D_yu\Big).$$
The above relations shows that $\tilde u\in W^{2,p}_{v}(\alpha, \alpha,m)$ if and only if  $u\in W^{2,p}_{\mathcal N}(\alpha, \alpha,m)$. This proves the required claim. The last claim follows by the density of $\mathcal{C}$ in $W^{2,p}_{\mathcal N}(\alpha, \alpha,m)$ since $\mathcal{C}_v=	S_{0,-\frac d c}\left(\mathcal{C}\right)$.\\
\qed

As in Propositions \ref{Trace D_yu in W} and \ref{trace u in W op} we can  provide  an equivalent description of $W^{2,p}_{v}(\alpha, \alpha, m)$, adapted to the degenerate  operator $D_{yy}+y^{-1}v\cdot \nabla $.  In the first formulation we shows that the  oblique boundary condition  in the integral  form $y^{\alpha-1}v\cdot \nabla u\in L^p_m$ is actually equivalent to the trace condition $\ds \lim_{y\to 0}y^c v\cdot \nabla u=0$. 
\begin{prop}\label{Equiv Wv alpha eq}
	Let    $\frac{m+1}{p}<c+1-\alpha$.  Then  one has 
	\begin{align*} 
		W^{2,p}_{v}(\alpha, \alpha, m)=&\left\{u \in  W^{2,p}_{loc}(\R^{N+1}_+): u,\ y^{\alpha}\Delta_xu\in L^p_m  \right. \\[1ex]			
		&\left.\hspace{10ex} y^{\alpha}D_{yy}u+y^{\alpha-1}v\cdot \nabla u \in L^p_m\text{\;\;and\;\;}\lim_{y\to 0}y^c v\cdot \nabla u=0\right\}
	\end{align*}
	and the norms $\|u\|_{W^{2,p}_{v}(\alpha,\alpha,m)}$ and $$\|u\|_{L^p_m}+\|y^{\alpha}\Delta_x u\|_{L^p_m}+\|y^{\alpha}D_{yy}u+y^{\alpha-1}v\cdot \nabla u\|_{L^p_m}$$ are equivalent on $W^{2,p}_{ v}(\alpha_1, \alpha_2, m)$.
	Finally,  when $0<\frac{m+1}p\leq c-1$ then 
	\begin{align*} 
		W^{2,p}_{v}(\alpha, \alpha, m)=&\left\{u \in  W^{2,p}_{loc}(\R^{N+1}_+): u,\;\; y^{\alpha}\Delta_xu,\;\;   y^{\alpha}D_{yy}u+y^{\alpha-1}v\cdot \nabla u \in L^p_m\right\}.
	\end{align*}
\end{prop}
{\sc{Proof.}}
Let $\tilde u\in W^{2,p}_{v}(\alpha, \alpha,m)$. By Proposition \ref{Equi obli Neumm sob alpha}, $\tilde u={S_{0,-\frac d c}}u$  where   $u\in W^{2,p}_{\mathcal N}(\alpha, \alpha,m)$. We observe that by Proposition \ref{Trace D_yu in W}, $u$ is characterized by 
\begin{align}\label{Equiv Wv alpha eq 1}
	u,\;\; y^{\alpha}\Delta_xu,\quad y^{\alpha}D_{yy}u+cy^{\alpha-1}D_yu \in L^p_m\quad \text{and\quad }\lim_{y\to 0}y^c D_yu=0.
\end{align}
We remark preliminarily that the Calderon-Zygmund inequality (see e.g. \cite[Theorem 2, Chapter 4, Section 1]{krylov-Lp}) yields
\begin{align*}
	\int_{\R^N} |D_{x_i x_j}u(x,y)|^p\,dx\leq C\int_{\R^N} |\Delta_xu(x,y)|^p\,dx.
\end{align*}
Multiplying by $y^{p\alpha+m
}$ and integrating over $\R_+$ 
we obtain 
\begin{align}\label{Equiv Wv alpha eq 2}
	\sum_{i,j=1}^n\|y^{\alpha} D_{x_ix_j}u\|_{L^p_m}\leq C\|y^{\alpha} \Delta_x u\|_{L^p_m}
\end{align}
and the same relation obviously holds for $\tilde u$. Moreover by Proposition \ref{Trace D_yu in W} again we also obtain
\begin{align}\label{Equiv Wv alpha eq 3}
	\|y^{\alpha}\nabla_xD_yu\|_{L^p_m}\leq C\Big[\|u\|_{L^p_m}+\|y^{\alpha}\Delta_x u\|_{L^p_m}+\|y^{\alpha}D_{yy}u+cy^{\alpha-1}D_yu\|_{L^p_m}\Big].
\end{align}
By Proposition \ref{Isometry shift} and Corollary \ref{Isometry action 2 spe}  with $\beta=0$ and $\omega=-\frac d c$ we have
\begin{itemize}
	\item [1.] $y^{\alpha}D_{x_ix_j}\tilde u={S_{0,-\frac d c}} \left(y^{\alpha}D_{x_ix_j} u\right)$;\medskip
	\item[2.]   $y^{\alpha-1} v \cdot \nabla \tilde u={S_{0,-\frac d c}}\Big(cy^{\alpha-1} D_yu\Big)$;\medskip
	\item[3.]	$y^{\alpha}D_{yy} \tilde u+y^{\alpha-1}v\cdot\nabla \tilde u={S_{0,-\frac d c}}\Big(\frac{1}{c^2}y^{\alpha}\left(D^2_xu\cdot d,d\right)-\frac{2}cy^{\alpha}\left(d,\nabla_xD_yu\right)+y^{\alpha}D_{yy}u+cy^{\alpha-1}D_yu\Big)$.
\end{itemize}
The above relations and  \eqref{Equiv Wv alpha eq 2}, \eqref{Equiv Wv alpha eq 3}, then  shows that the requirements in \eqref{Equiv Wv alpha eq 1} are equivalent to 
\begin{align*}
	\tilde u,\quad y^{\alpha}\Delta_x\tilde u,\;\; y^{\alpha}D_{yy}\tilde u+y^{\alpha-1}v\cdot \nabla \tilde u \in L^p_m\text{\quad and\quad}\lim_{y\to 0}y^c v\cdot \nabla \tilde  u=0
\end{align*}
which proves the first required claim. The claim for $0<\frac{m+1}p\leq c-1$ follows similarly.
\qed

The following proposition involves a Dirichlet, rather than an oblique, boundary condition,  in a certain range of parameters.
\begin{prop}\label{Equiv Wv alpha dir}
	Let   $c\geq 1$ and $\frac{m+1}{p}<c+1-\alpha$.  The following properties hold.
	\begin{itemize}
		\item[(i)] If $c>1$ then 
	\begin{align*} 
	W^{2,p}_{v}(\alpha, \alpha, m)=&\left\{u \in  W^{2,p}_{loc}(\R^{N+1}_+): u,\ y^{\alpha}\Delta_xu\in L^p_m  \right. \\[1ex]			
	&\left.\hspace{10ex} y^{\alpha}D_{yy}u+y^{\alpha-1}v\cdot \nabla u \in L^p_m\text{\;\;and\;\;}\lim_{y\to 0}y^{c-1} u=0\right\}
\end{align*}
		\item[(ii)] If $c=1$ then 
		\begin{align*} 
			W^{2,p}_{v}(\alpha, \alpha, m)=&\left\{u \in  W^{2,p}_{loc}(\R^{N+1}_+): u,\ y^{\alpha}\Delta_xu\in L^p_m  \right. \\[1ex]			
			&\left.\hspace{10ex} y^{\alpha}D_{yy}u+y^{\alpha-1}v\cdot \nabla u \in L^p_m\text{\;\;and\;\;}\lim_{y\to 0}y^{c-1} u\in\C\right\}
		\end{align*}
	\end{itemize} 
\end{prop}
{\sc{Proof.}}
The proof follows as in Proposition \ref{Equiv Wv alpha eq} by using  Proposition \ref{trace u in W op} in place of Proposition \ref{Equi obli Neumm sob alpha}.
\qed 
\subsubsection{The general case: the space $W^{2,p}_{v}(\alpha_1,\alpha_2,m)$} \label{sobolev oblique alpha diff section}
We now extend the previous definition also to the case $\alpha_1\neq\alpha_2$. In contrast to the case $\alpha_1=\alpha_2$, a distortion correction defined by the coefficient
\begin{align*}
	\beta_{\alpha}\overset{def}=\frac{\alpha_1-\alpha_2}{2}
\end{align*} 
appears in the definition of $W^{2,p}_{v}(\alpha_1,\alpha_2,m)$.   This  is due to the possibly  different weights $y^{\alpha_1}$, $y^{\alpha_2}$ which appear in the $x$ and $y$ directions. It does not appears, obviously,  in the case $\alpha_1=\alpha_2$, where $\beta_\alpha=0$, but also when $d=0$ (see Remark \ref{oss defi W_v}).  This correction is  essential for the validity of equalities \eqref{Tb, Sb eq 1} and \eqref{Tb, Sb eq 2} and for characterizing,  in Section \ref{Consequenses}, the domain of the degenerate operator $\mathcal L$ defined in \eqref{def operator sobolev section}.

 We start by defining
\begin{align*}
	F^{2,p}(\alpha_1, \alpha_2,m)=&\left\{u\in W^{2,p}_{loc}(\R^{N+1}_+):\ u,\  y^{\alpha_1} D_{x_ix_j}u,\ y^\frac{\alpha_1}{2} D_{x_i}u, \,y^\frac{\alpha_1+\alpha_2}{2} D_{y}\nabla_x u,  \right.\\[1ex]
	&\left.\hspace{24ex}y^{\alpha_2}\left(D_{yy}u+\beta_{\alpha}\frac{D_yu}{y}\right),\ y^{\frac{\alpha_2}{2}}D_{y}u\in L^p_m\right\}
\end{align*}
with the norm
\begin{align*}
	\|u\|_{F^{2,p}(\alpha_1,\alpha_2,m)}{=}\|u\|_{L^p_m}&+\sum_{i,j=1}^n\|y^{\alpha_1} D_{x_ix_j}u\|_{L^p_m}+\sum_{i=1}^n\|y^{\frac{\alpha_1}2} D_{x_i}u\|_{L^p_m}+\|y^\frac{\alpha_1+\alpha_2}{2} D_{y}\nabla_x u\|_{L^p_m}\\[1ex]
	&+\left\|y^{\alpha_2}\left(D_{yy}u+\beta_{\alpha}\frac{D_yu}y\right)\right\|_{L^p_m}+\|y^{\frac{\alpha_2}{2}}D_{y}u\|_{L^p_m}.
\end{align*}
 Next we impose   the weighted  oblique derivative boundary condition  
 $$\left(y^{\beta_\alpha}\, d \cdot \nabla_x u+cD_yu\right)(x,0)=0$$
  (see Proposition \ref{Equiv Wv alpha diff}) in the integral form 
\begin{align}
	y^{\frac{\alpha_1+\alpha_2}2-1} d \cdot \nabla_x u+cy^{\alpha_2-1}D_yu=	y^{\alpha_2-1}\Big( y^{\beta_\alpha}\,d \cdot \nabla_x u+cy^{\alpha_2-1}D_yu\Big)\in L^p_m.
\end{align}
Accordingly we define
\begin{align*}
	W^{2,p}_{v}(\alpha_1,\alpha_2,m)\stackrel{def}{=}=\{u \in F^{2,p}(\alpha_1, \alpha_2,m):\  y^{\frac{\alpha_1+\alpha_2}2-1} d \cdot \nabla_x u+cy^{\alpha_2-1}D_yu \in L^p_m\}
\end{align*}
with the norm
$$
\|u\|_{W^{2,p}_{v}(\alpha_1,\alpha_2,m)}\stackrel{def}{=}\|u\|_{F^{2,p}(\alpha_1,\alpha_2,m)}+\|y^{\frac{\alpha_1+\alpha_2}2-1} d \cdot \nabla_x u+cy^{\alpha_2-1}D_yu\|_{ L^p_m}.
$$
\begin{os}\label{oss defi W_v}\mbox{}
	\begin{itemize}
		\item[(i)] The difference between the spaces 	$W^{2,p}(\alpha_1,\alpha_2,m)$ and $	F^{2,p}(\alpha_1,\alpha_2,m)$ relies on the requirement  $y^{\alpha_2}\left(D_{yy}u+\beta_\alpha\frac{D_yu}{y}\right)\in\L^p_m$ which, in the definitions of $F^{2,p}(\alpha_1,\alpha_2,m)$ and of $W^{2,p}_{v}(\alpha_1,\alpha_2,m)$, cannot be split into 
		\begin{align*}
			y^{\alpha_2}D_{yy}u,\; y^{\alpha_2-1}D_yu\in\L^p_m
		\end{align*}
		as in the definition of $W^{2,p}_{\mathcal N}(\alpha_1,\alpha_2,m)$.
		\item[(ii)] Both the requirements 
		\begin{align}\label{eq BC obliq diff}
			y^{\alpha_2}\left(D_{yy}u+\beta_\alpha\frac{D_yu}{y}\right)\in\L^p_m,\qquad y^{\frac{\alpha_1+\alpha_2}2-1} d \cdot \nabla_x u+cy^{\alpha_2-1}D_yu\in L^p_m
		\end{align}
		in the definition of $W^{2,p}_{v}(\alpha_1,\alpha_2,m)$ are essential for the validity of Propositions \ref{Equi obli Neumm sob} and  \ref{Isometry action sob} (see Remark \ref{oss necessity Tb}).
		\item[(iii)] The hypotheses  $y^{\alpha_2}\left(D_{yy}u+\beta_\alpha\frac{D_yu}{y}\right)\in\L^p_m$ is equivalent to 
		\begin{align*}
			y^{\alpha_2}D_{yy}u-\frac{\beta_\alpha}{c}y^{\frac{\alpha_1+\alpha_2}2-1} d \cdot \nabla_x u\in\L^p_m.
		\end{align*}
		This follows by combining linearly the two conditions in \eqref{eq BC obliq diff}. In particular, when $d=0$, one has in any case $$W^{2,p}_{(0,c)}(\alpha_1,\alpha_2,m)=W^{2,p}_{\mathcal N}(\alpha_1,\alpha_2,m),\qquad \text{for any}\quad c\in\R.$$
	\item[(iv)] The boundary condition  \eqref{eq bound cond alpha diff} is  defined up to a normalization constant i.e.
	\begin{align*}
		W^{2,p}_{\mu v}(\alpha_1,\alpha_2,m)=W^{2,p}_{v}(\alpha_1,\alpha_2,m),\qquad \forall \mu\in\R\setminus\{0\}.
	\end{align*}
	\end{itemize}
\end{os}
\medskip
Note that in the special case $\alpha_1=\alpha_2$ the above definition is consistent with the one given in Subsection \ref{section Wv alpha eq}. Moreover  Remark \ref{oss defi W_v} (iii) justifies the following definition.
\begin{defi}\label{convention Wv Wn alpha}
	  To shorten some statements we also write
\begin{align*}
	W^{2,p}_{(0,0)}(\alpha_1,\alpha_2,m)\stackrel{def}{=}W^{2,p}_{\mathcal N}(\alpha_1,\alpha_2,m).
\end{align*}
With this notation $W^{2,p}_{v}(\alpha_1,\alpha_2,m)$ is well defined  for any $v=(d,c)\in\R^{N+1}$ such that $d=0$ if $c=0$.  Note that  \begin{align*}
	W^{2,p}_{v}(\alpha_1,\alpha_2,m)=W^{2,p}_{\mathcal N}(\alpha_1,\alpha_2,m),\qquad \text{when\quad $c=0$\quad or \quad  $d=0$}.
\end{align*}  
Under this identification, we also define
\begin{align} \label{defC oblique}
 \mathcal{C}_v
	&:=\left \{u \in C_c^\infty \left(\R^N\times[0, \infty)\right), \ \left(y^{\beta_\alpha} d \cdot \nabla_x u+cD_yu\right)(x,y)=0\  {\rm for} \ y \leq \delta\ {\rm  and \ some\ } \delta>0\right \}
\end{align}
with the convention that $\mathcal{C}_{v}:=\mathcal C$ when $v=(0,0)$ and 	$\mathcal{C}$ is the set defined in \eqref{defC}. 

\end{defi}

\medskip

As before we clarify the relation  between the spaces $W^{2,p}_{N}(\alpha_1,\alpha_1,m)$ and   $W^{2,p}_{v}(\alpha_1,\alpha_2,m)$ by generalizing Proposition \ref{Equi obli Neumm sob alpha}. $W^{2,p}_{v}(\alpha_1, \alpha_2,m)$ is related to $W^{2,p}_{\mathcal N}(\alpha_1,\alpha_2,m)$ by means of the isometry $S_{\beta_\alpha,\omega}$ of $L^p_m$   defined in \eqref{Tran def} with
 \begin{align*}
	\beta_\alpha=\frac{\alpha_1-\alpha_2}2\neq -1, \qquad \omega=-\frac d {c(\beta_\alpha+1)}, 
	\end{align*}  
 namely
\begin{align*}
	S_{\beta_\alpha,\omega}\, u(x,y)&:=u\left(x-\omega y^{\beta_\alpha+1},y\right),\quad (x,y)\in\R_+^{N+1}.
\end{align*}

\begin{prop}\label{Equi obli Neumm sob}
	Let $\beta_\alpha=\frac{\alpha_1-\alpha_2}2\neq -1$ and  $\omega=-\frac d {c(\beta_\alpha+1)}$. Then one has 	
	\begin{align}\label{Tb, Sb eq 1}
		S_{\beta_\alpha,\omega}\,\left(W^{2,p}_{\mathcal N}(\alpha_1, \alpha_2,m)\right)	=W^{2,p}_{v}(\alpha_1, \alpha_2,m)
	\end{align}
	In particular the set $\mathcal{C}_v$ defined in \eqref{defC oblique} 	is dense in $W^{2,p}_{v}(\alpha_1, \alpha_2,m)$.
\end{prop}
{\sc{Proof.}} Let $u\in W^{2,1}_{loc}\left(\R^{N+1}_+\right)$  and let us set $\tilde u={S_{\beta_\alpha,\omega}}u$. Then  by Proposition \ref{Isometry shift} one has 
\begin{itemize}
	\item [1.] $y^{\alpha_1}D_{x_ix_j}\tilde u={S_{\beta_\alpha,\omega}} \left(y^{\alpha_2}D_{x_ix_j} u\right)$, \quad $y^{\frac{\alpha_1}2}D_{x_i}\tilde u={S_{\beta_\alpha,\omega}}\left(y^{\frac{\alpha_1}2}D_{x_i} u\right)$;\medskip
	\item[2.]  $y^{\frac{\alpha_2}2}D_y \tilde u={S_{\beta_\alpha,\omega}}\Big((\beta_\alpha+1)y^{\frac{\alpha_1}2} (\nabla_x u,\omega)+y^{\frac{\alpha_2}2}D_yu\Big)$,	\\[1ex]
	$y^{\alpha_2-1}D_y \tilde u={S_{\beta_\alpha,\omega}}\Big((\beta_\alpha+1)y^{\frac{\alpha_1+\alpha_2}2-1} (\nabla_x u,\omega)+y^{\alpha_2-1}D_yu\Big)$;\medskip
	\item[3.]	 $y^{\alpha_2}D_{yy} \tilde u={S_{\beta_\alpha,\omega}}\Big((\beta_\alpha+1)^2y^{\alpha_1} (D^2_x u\cdot\omega,\omega)+2(\beta_\alpha+1)y^{\frac{\alpha_1+\alpha_2}2} (\nabla_x D_y u, \omega)$
	\begin{flushright}	\vspace{-2ex}$+\beta_\alpha (\beta_\alpha+1)y^{\frac{\alpha_1+\alpha_2}2-1}(\nabla_x u,\omega)+ y^{\alpha_2}D_{yy}u\Big)$;\end{flushright}
	\item[4.] $y^{\frac{\alpha_1+\alpha_2}2}\nabla_xD_{y} \tilde u={S_{\beta_\alpha,\omega}}\Big((\beta_\alpha+1)y^{\alpha_1} D^2_x u\cdot \omega+y^{\frac{\alpha_1+\alpha_2}2}\nabla_xD_yu\Big)$;
	\medskip
	\item[5.]$y^{\frac{\alpha_1+\alpha_2}2-1} d \cdot \nabla_x \tilde u+cy^{\alpha_2-1}D_y\tilde u ={S_{\beta_\alpha,\omega}}\Big(cy^{\alpha_2-1} D_yu\Big)$.
\end{itemize}
The above relations shows that $\tilde u\in W^{2,p}_{v}(\alpha_1, \alpha_2,m)$ if and only if  $u\in W^{2,p}_{\mathcal N}(\alpha_1, \alpha_2,m)$. This proves the required claim. The last claim follows by the density of $\mathcal{C}$ in $W^{2,p}_{\mathcal N}(\alpha_1, \alpha_2,m)$ since $\mathcal{C}_v=	S_{\beta,\omega}\left(\mathcal{C}\right)$.\\
\qed

\begin{os}
	Recalling Remark \ref{oss defi W_v} (iv), equality \eqref{Tb, Sb eq 1} can be written equivalently as	\begin{align}\label{Tb, Sb eq 1bis}
		S_{\beta_\alpha,\omega}\,\left(W^{2,p}_{\mathcal N}(\alpha_1, \alpha_2,m)\right)	=W^{2,p}_{v}(\alpha_1, \alpha_2,m), \qquad v=\left(-\omega(\beta_\alpha+1),1\right)
	\end{align}
which is valid for any $\omega\in\R^N$.
\end{os}
As in Propositions \ref{Equiv Wv alpha eq} and  \ref{Equiv Wv alpha dir} we  provide  an equivalent description of $W^{2,p}_{v}(\alpha_1, \alpha_2, m)$, adapted to the degenerate  operator 
\begin{align}\label{degenerate B oblique}
	\mathcal A:&=y^{\alpha_2}D_{yy}+y^{\frac{\alpha_1+\alpha_2}2-1} \tilde d \cdot \nabla_x + \mu y^{\alpha_2-1}D_y,\qquad \tilde d=\frac{\mu+\beta_\alpha}{c}d.
\end{align}

\begin{prop}\label{Equiv Wv alpha diff}
	Let $\mu\in\R$ such that   $\frac{m+1}{p}<\mu+1-\alpha_2$ and   $\mathcal A$ be the operator in \eqref{degenerate B oblique}.
	 Then  one has 
	\begin{align*} 
		W^{2,p}_{v}(\alpha_1, \alpha_2, m)=&\left\{u \in  W^{2,p}_{loc}(\R^{N+1}_+): u,\ y^{\alpha_1}\Delta_xu,\; \mathcal Au \in L^p_m \right. \\[1ex]			
		&\left.\hspace{17ex}\text{\;\;and\;\;}\lim_{y\to 0}y^\mu \left(y^{\beta_\alpha}\, d \cdot \nabla_x u+cD_yu\right)=0\right\}
	\end{align*}
	and the norms $\|u\|_{W^{2,p}_{v}(\alpha_1,\alpha_2,m)}$ and $$\|u\|_{L^p_m}+\|y^{\alpha_1}\Delta_x u\|_{L^p_m}+\|	\mathcal Au\|_{L^p_m}$$ are equivalent on $W^{2,p}_{v}(\alpha_1, \alpha_2, m)$.
	Finally,  when $0<\frac{m+1}p\leq \mu-1$ then 
	\begin{align*} 
		W^{2,p}_{v}(\alpha_1, \alpha_2, m)=&\left\{u \in  W^{2,p}_{loc}(\R^{N+1}_+): u,\ y^{\alpha_1}\Delta_xu,\; 	\mathcal A u \in L^p_m\right\}.
	\end{align*}
\end{prop}
{\sc{Proof.}} Let us suppose, preliminarily, $\mu=c$. The proof, in this case, follows as in  Proposition \ref{Equiv Wv alpha eq} by using the isometry $S_{\beta_\alpha,\omega}$ in place of $S_{0,-\frac d c}$. The claim for a general $\mu\in\R$ follows, recalling Remark \ref{oss defi W_v} (iv),  by writing 
\begin{align*}
		W^{2,p}_{v}(\alpha_1, \alpha_2, m)=	W^{2,p}_{(\frac\mu c d,\mu)}(\alpha_1, \alpha_2, m)
\end{align*}
and by using the previous step with $(\frac\mu c d,\mu)$ in place of v .\qed

\begin{prop}\label{Equiv Wv alpha diff dir}
	Let   $\mu\geq 1$, $\frac{m+1}{p}<\mu+1-\alpha_2$  and   $\mathcal A$ be the operator in \eqref{degenerate B oblique}.   The following properties hold.
	\begin{itemize}
		\item[(i)] If $\mu>1$ then 
		\begin{align*} 
			W^{2,p}_{v}(\alpha_1, \alpha_2, m)=&\left\{u \in  W^{2,p}_{loc}(\R^{N+1}_+): u,\ y^{\alpha_1}\Delta_xu,\   	\mathcal Au \in L^p_m\text{\;\;and\;\;}\lim_{y\to 0}y^{\mu-1} u=0\right\}
		\end{align*}
		\item[(ii)] If $\mu =1$ then 
		\begin{align*} 
			W^{2,p}_{v}(\alpha, \alpha, m)=&\left\{u \in  W^{2,p}_{loc}(\R^{N+1}_+): u,\ y^{\alpha_1}\Delta_xu,\ 	\mathcal A u \in L^p_m\text{\;\;and\;\;}\lim_{y\to 0}y^{\mu-1} u\in\C\right\}
		\end{align*}.
	\end{itemize} 
\end{prop}
 {\sc{Proof.}} The proof follows as in Propositions  \ref{Equiv Wv alpha dir} and \ref{Equiv Wv alpha diff}. \qed

\medskip

The next Proposition shows that in a certain range of parameters the spaces $W^{2,p}_{\mathcal N}(\alpha_1, \alpha_2, m)$ and  $W^{2,p}_{v}(\alpha_1, \alpha_2, m)$ coincide.

\begin{prop}\label{Hardy Rellich Sob}
	If $\frac{m+1}p>1-\frac{\alpha_1+\alpha_2}{2}$ then $W^{2,p}_{\mathcal N}(\alpha_1, \alpha_2, m)= W^{2,p}_{v}(\alpha_1, \alpha_2, m)$ and one has 
	\begin{align*}
		\|y^{\frac{\alpha_1+\alpha_2}{2}-1}\nabla_{x}u\|_{L^p_m}\leq C \|y^\frac{\alpha_1+\alpha_2}{2} D_{y}\nabla_x u\|_{L^p_m},\qquad \forall u\in  W^{2,p}_{\mathcal N}(\alpha_1, \alpha_2, m).
	\end{align*} 
\end{prop}
{\sc{Proof.}}  Let  $u\in C_c^\infty \left(\R^N\times[0, \infty)\right)$. By applying \cite[Proposition 3.2 (iii)]{MNS-Sobolev} to $\nabla_x u$ we get 
\begin{align*}
	\|y^{\frac{\alpha_1+\alpha_2}{2}-1}\nabla_{x}u\|_{L^p_m}\leq C \|y^\frac{\alpha_1+\alpha_2}{2} D_{y}\nabla_x u\|_{L^p_m}.
\end{align*} 
Since by Theorem \ref{core gen} and Proposition \ref{Equi obli Neumm sob}, the set $ C_c^\infty \left(\R^N\times[0, \infty)\right)$ is dense in both $W^{2,p}_{\mathcal N}(\alpha_1, \alpha_2, m)$ and  $W^{2,p}_{v}(\alpha_1, \alpha_2, m)$, the above inequality extends to both spaces. This  shows in particular that for function $u$ belonging to  $W^{2,p}_{\mathcal N}(\alpha_1, \alpha_2, m)$ or $W^{2,p}_{v}(\alpha_1, \alpha_2, m)$, so that in particular $y^\frac{\alpha_1+\alpha_2}{2} D_{y}\nabla_x u\in L^p_m$, one has 
\begin{align*}
	y^{\frac{\alpha_1+\alpha_2}2-1} d \cdot \nabla_x u+cy^{\alpha_2-1}D_yu \in L^p_m\qquad \Longleftrightarrow\qquad y^{\alpha_2-1}D_yu \in L^p_m.
\end{align*}
This, recalling (iii) in Remark \ref{oss defi W_v},  proves the required claim.
\qed 
\medskip

We end the section by showing  how the spaces  $W^{2,p}_{\mathcal{N}}(\alpha_1,\alpha_2,m)$, $W^{2,p}_{v}(\alpha_1,\alpha_2,m)$, introduced so far, transform under the action of the  map   \eqref{Gen Kelvin def} with $k=0$,  $\beta\neq -1$, namely
\begin{align}
	\label{Gen Kelvin def k=0}
	T_{0,\beta}u(x,y)&:=|\beta+1|^{\frac 1 p}u(x,y^{\beta+1}),\quad (x,y)\in\R^{N+1}_+.
\end{align}
Observe that
$$ T_{0,\beta}^{-1}=T_{0,-\frac{\beta}{\beta+1}\,}.$$

\begin{prop}\label{Isometry action sob} Let $1\leq p\leq \infty$, $\beta \in\R$, $\beta\neq -1$ and $m\in\R$. The following properties hold.
	\begin{itemize}
		\item[(i)] $T_{0,\beta}$ maps isometrically  $L^p_{\tilde m}$ onto $L^p_m$  where $ \tilde m=\frac{m-\beta}{\beta+1}$.
		\item[(ii)] Setting $\tilde\alpha_1=\frac{\alpha_1}{\beta+1}$,\quad $\tilde\alpha_2=\frac{\alpha_2+2\beta}{\beta+1}$ one has  
		\begin{align}\label{Tb, Sb eq 2}
			\nonumber W^{2,p}_{\mathcal{N}}(\alpha_1,\alpha_2,m)&=T_{0,\beta}\left(W^{2,p}_{\mathcal{N}}(\tilde \alpha_1,\tilde\alpha_2,\tilde m)\right),\\[1.5ex]
			W^{2,p}_{v}(\alpha_1,\alpha_2,m)&=T_{0,\beta}\left(W^{2,p}_{\tilde v}(\tilde \alpha_1,\tilde\alpha_2,\tilde m)\right),\qquad \tilde v=\big(d, c(\beta+1)\big).
		\end{align}
	\end{itemize}
	In particular choosing $\beta=\beta_\alpha=\frac{\alpha_1-\alpha_2}2$ and setting $\tilde \alpha=\frac{2\alpha_1}{\alpha_1-\alpha_2+2}$ and $m=\frac{2m-\alpha_1+\alpha_2}{\alpha_1-\alpha_2+2}$ one has
	\begin{align*}
		W^{2,p}_{v}(\alpha_1,\alpha_2,m)&=T_{0,\beta_\alpha}\left(W^{2,p}_{\tilde v}(\tilde \alpha,\tilde \alpha ,\tilde m)\right),\qquad \tilde v=\big(d, c(\beta_\alpha+1)\big).
	\end{align*}
\end{prop}
{\sc{Proof.}}  (i) and the first equality in (ii)  follow by a straightforward application  of Proposition \ref{Isometry action der} with $k=0$ as in  \cite[Proposition 2.2]{MNS-Sobolev}. 
To prove the second equality in (ii) we set $\tilde v=(d,\tilde c):=(d,c(\beta+1))$ and 
\begin{align*}
	\gamma=\frac{\alpha_1-\alpha_2}2,\qquad \omega=-\frac d {c(\gamma+1)},\qquad \tilde \gamma=\frac{\tilde\alpha_1-\tilde\alpha_2}2,\qquad \tilde\omega=-\frac d {\tilde c(\tilde\gamma+1)}.
\end{align*}
and we observe preliminarily that, by construction, one has 
\begin{align*}
	\tilde \gamma+1=\frac{\gamma+1}{\beta+1},\qquad \tilde \omega=\omega.
\end{align*}
Then, recalling the definitions of  	$S_{\gamma,\omega}$ and  $T_{0,\beta}$,   a straightforward calculation yields
\begin{align*}
	S_{\gamma,\omega}\circ T_{\beta}=T_{\beta} \circ S_{\tilde\gamma,\tilde\omega}.
\end{align*}
Indeed the latter equality follows by the previous relations on the coefficients after observing that
\begin{align*}
S_{\gamma,\omega}\left(T_{\beta} u\right)(x,y)=\left(T_{\beta} u\right)\left(x+\omega y^{\gamma+1},y\right)=|\beta+1|^{\frac 1 p}u\left(x+\omega y^{\gamma+1},y^{\beta+1}\right)
\end{align*}
and
\begin{align*}
	T_{\beta}\left(S_{\tilde\gamma,\tilde \omega} u\right)(x,y)=|\beta+1|^{\frac 1 p}(S_{\tilde\gamma,\tilde \omega}u)\left(x,y^{\beta+1}\right)=|\beta+1|^{\frac 1 p}u\left(x+\tilde\omega y^{(\tilde \gamma+1)(\beta+1)},y^{\beta+1}\right).
\end{align*}
Then using the first equality in (ii) and Proposition \ref{Equi obli Neumm sob}, one has 	
\begin{align*}
	W^{2,p}_{v}(\alpha_1,\alpha_2,m)&=S_{\gamma,\omega}\,\left(W^{2,p}_{\mathcal N}(\alpha_1, \alpha_2,m)\right)=S_{\gamma,\omega}\,\left(T_{0,\beta}\left(W^{2,p}_{\mathcal{N}}(\tilde \alpha_1,\tilde\alpha_2,\tilde m)\right)\right) \\[1.5ex]
	&=T_{0,\beta}\left(S_{\tilde\gamma,\tilde\omega}\,\left(W^{2,p}_{\mathcal N}(\tilde\alpha_1, \tilde\alpha_2,m)\right)\right)	=T_{0,\beta}\left(W^{2,p}_{\tilde v}(\tilde \alpha_1,\tilde\alpha_2,\tilde m)\right).
\end{align*}
This proved the required claim.
\qed

\begin{os}\label{oss necessity Tb}
	It is essential to deal with  $W^{2,p}_{\mathcal{N}}(\alpha_1,\alpha_2,m)$, $W^{2,p}_{v}(\alpha_1,\alpha_2,m)$: in  general the map $T_{0,\beta}$ does  not  transform  $W^{2,p}(\tilde \alpha_1,\tilde\alpha_2,\tilde m)$ into $W^{2,p}(\alpha_1,\alpha_2,m)$ since by (ii)-3 of Proposition  \ref{Isometry action der} one has 
	$$y^{\alpha_2}D_{yy} (T_{\beta\,} u)=(\beta+1)T_{\beta\,}\Big((\beta+1)y^{{\tilde\alpha}_2}D_{yy}u+\beta y^{{\tilde\alpha}_2-1}D_y u\Big).$$
\end{os}

\subsection{The space $W^{2,p}_{\mathcal R}(\alpha_1,\alpha_2,m)$}\label{section W_r}

We  consider also an integral version of the Dirichlet boundary condition, namely  a weighted summability requirement for $y^{-2}u$ and introduce 
$$
W^{2,p}_{\mathcal R}(\alpha_1, \alpha_2, m)\stackrel{def}{=}\{u \in  W^{2,p}(\alpha_1, \alpha_2, m): y^{\alpha_2-2}u \in L^p_m\}
$$
with the norm $$\|u\|_{W^{2,p}_{\mathcal R}(\alpha_1, \alpha_2, m)}\stackrel{def}{=}\|u\|_{W^{2,p}(\alpha_1, \alpha_2, m)}+\|y^{\alpha_2-2}u\|_{L^p_m}.$$
The symbol $\mathcal R$ stands for "Rellich", since Rellich inequalities concern with the summability of $y^{-2}u$.

\begin{prop} \label{RN} The following properties hold.
	\begin{itemize}
		\item[{\rm(i)}]  For any $u \in W^{2,p}_{\mathcal R}(\alpha_1, \alpha_2, m)$ one has 
		\begin{align*}
			\|y^{\alpha_2-1}D_yu\|_{L^p_m}&\leq C\left( \|y^{\alpha_2}D_{yy}u\|_{L^p_m}+\|y^{\alpha_2-2}u\|_{L^p_m}\right),\\[1.5ex]
			\|y^{\frac{\alpha_1+\alpha_2}2-1}\nabla_x u\|_{L^p_m}&\leq C\left( \|y^{\alpha_1}D_{xx}u\|_{L^p_m}+\|y^{\alpha_2-2}u\|_{L^p_m}\right)  
		\end{align*}
		In particular $W^{2,p}_{\mathcal R}(\alpha_1, \alpha_2, m)\subseteq W^{2,p}_{\mathcal N}(\alpha_1, \alpha_2, m)\cap W^{2,p}_{v}(\alpha_1, \alpha_2, m)$.
		\item[{\rm(ii)}]   $C_c^\infty (\R^{N+1}_+)$ is dense in $W^{2,p}_{\mathcal R}(\alpha_1, \alpha_2, m)$.
		\item[{\rm(iii)}] If $\frac{m+1}{p}>2-\alpha_2$, then 
		$$W^{2,p}_{\mathcal R}(\alpha_1, \alpha_2, m) =  W^{2,p}_{\mathcal N}(\alpha_1, \alpha_2, m)=  W^{2,p}_{v}(\alpha_1, \alpha_2, m)=W^{2,p}(\alpha_1, \alpha_2, m)$$ with equivalence of the corresponding norms. 
	\end{itemize}
\end{prop}
{\sc{Proof.}} The proof follows by \cite[Proposition 5.1]{MNS-Sobolev}. We need only to prove the second inequality in (i). With this aim let $u \in W^{2,p}_{\mathcal R}(\alpha_1, \alpha_2, m)$. We	use the classic interpolative inequality
$$\|\nabla_x u(\cdot,y)\|_{L^p(\R^N)}\leq \varepsilon \|\Delta_x u(\cdot,y)\|_{L^p(\R^N)}+\frac {C(N,p)} \varepsilon \| u(\cdot,y)\|_{L^p(\R^N)}.$$
Multiplying the above inequality by $y^{\frac{\alpha_1+\alpha_2}2-1}$ and choosing  $\varepsilon=y^{1-\frac{\alpha_2-\alpha_1}2}$ we  get 
$$\|y^{\frac{\alpha_1+\alpha_2}2-1}\nabla_x u(\cdot,y)\|_{L^p(\R^N)}\leq  \|y^{\alpha_1}\Delta_x u(\cdot,y)\|_{L^p(\R^N)}+C(N,p)  \|y^{\alpha_2-2} u(\cdot,y)\|_{L^p(\R^N)}.$$
The required estimate then follows after raising to the power $p$ and integrating in $y$. We remark that in (iii), under the range of parameters \eqref{param sobolev}, we have $\frac{m+1}{p}>2-\alpha_2>1-\frac{\alpha_1+\alpha_2}2$ and then, by Proposition \ref{Hardy Rellich Sob},  $W^{2,p}_{\mathcal N}(\alpha_1, \alpha_2, m)=  W^{2,p}_{v}(\alpha_1, \alpha_2, m)$.\qed

Finally, we investigate the action of the isometry $T_{k,\beta}$ defined in \eqref{Gen Kelvin def}. We start with the case $k=0$.

\begin{prop}\label{Isometry action sob R}  Let $T_{0,\beta}$ the map  defined in \eqref{Gen Kelvin def k=0}. Then 
	\begin{align*}
		W^{2,p}_{\mathcal R}(\alpha_1,\alpha_2,m)=T_{0,\beta}\left(W^{2,p}_{\mathcal{R}}(\tilde \alpha_1,\tilde\alpha_2,\tilde m)\right)
	\end{align*}
	where  $ \tilde m=\frac{m-\beta}{\beta+1}$, $ \tilde\alpha_1=\frac{\alpha_1}{\beta+1}$, $\tilde\alpha_2=\frac{\alpha_2+2\beta}{\beta+1}$.  	In particular choosing $\beta=\beta_\alpha=\frac{\alpha_1-\alpha_2}2$ one has
	\begin{align*}
		W^{2,p}_{\mathcal R}(\alpha_1,\alpha_2,m)=T_{0,\beta_\alpha}\left(W^{2,p}_{\mathcal{R}}(\tilde \alpha,\tilde\alpha,\tilde m)\right),\quad m=\frac{2m-\alpha_1+\alpha_2}{\alpha_1-\alpha_2+2}, \quad \tilde \alpha=\frac{2\alpha_1}{\alpha_1-\alpha_2+2}.
	\end{align*}
\end{prop}
{\sc{Proof.}}
The claim follows from  Proposition \ref{Isometry action sob} since by Proposition \ref{RN} one has 
\begin{align*}
	W^{2,p}_{\mathcal R}(\alpha_1,\alpha_2,m)=W^{2,p}_{\mathcal N}(\alpha_1,\alpha_2,m)\cap \{u\in L^p_m: y^{\alpha_2-2}u\in L^p_m\} 
\end{align*}
and noticing that $y^{\alpha_2-2}u\in\ L^p_m$ if and only if $y^{{\tilde\alpha}_2-2}u\in\ L^p_{\tilde m}$.\qed

\medskip

We consider now the multiplication operator $T_{k,0}:u\mapsto y^ku$.

\begin{prop}{\em \cite[Lemma 5.3]{MNS-Sobolev}}\label{isometryRN}
	\label{y^k W}
	Let   $\alpha_2-\alpha_1<2$ and  $\frac{m+1}{p}>2-\alpha_2$. For every $k\in\R$,
	\begin{align*}
		T_{k,0}:  W^{2,p}_{\mathcal N}(\alpha_1, \alpha_2, m) \to  W^{2,p}_{\mathcal R}(\alpha_1, \alpha_2, m-kp)
	\end{align*}
	
	is an isomorphism {\rm (}we shall write $y^k  W^{2,p}_{\mathcal N}(\alpha_1, \alpha_2, m)= W^{2,p}_{\mathcal R}(\alpha_1, \alpha_2, m-kp)${\rm )}.
\end{prop}
\medskip

Finally we deal with the isometry of $L^p_m$,  $S_{\beta_\alpha,\omega}$ defined in \eqref{Tran def} with  $\beta_\alpha=\frac{\alpha_1-\alpha_2}2\neq -1$ and  $\omega=-\frac d {c(\beta_\alpha+1)}$.

\begin{prop}\label{Equi obli Diric sob}
	Let $\beta_\alpha=\frac{\alpha_1-\alpha_2}2\neq -1$ and  $\omega=-\frac d {c(\beta_\alpha+1)}$. Then one has 	
	\begin{align*}
		S_{\beta_\alpha,\omega}\,\left(W^{2,p}_{\mathcal R}(\alpha_1, \alpha_2,m)\right)	=W^{2,p}_{\mathcal R}(\alpha_1, \alpha_2,m)
	\end{align*}
\end{prop}
{\sc{Proof.}} The proof follows as in Proposition \ref{Isometry action sob R}  using Proposition \ref{Equi obli Neumm sob}.\qed

\section{The  operator with oblique boundary conditions} \label{Sect DOm obliqu}

In this section we study  parabolic problems related to the operator 
\begin{align*}
	\mathcal L&=y^{\alpha}\mbox{Tr }\left(AD^2\right)+y^{\alpha-1}\left(v,\nabla\right)
\end{align*}
defined in \eqref{def L transf alpha} in the  case  $b=0$. Here  $v=(d,c)\in\R^{N+1}$ and, under  the hypotheses in  Assumption \ref{assumption},  we always  assume $\alpha<2$ and the diffusion matrix 
$$A=	\left(
\begin{array}{c|c}
	Q  & { q}^t \\[1ex] \hline
	q& \gamma
\end{array}\right)$$
 to be symmetric and positive definite. We  endow $\mathcal L$  with Neumann and oblique boundary conditions in the sense specified below.  \\

We briefly  recall the weighted Sobolev spaces which are introduce and analysed in details, as well as their boundary conditions,  in Sections \ref{Section neumann sobolev} and \ref{section oblique sobolev}. 
We set, at first,
\begin{align*}
	W^{2,p}(\alpha,\alpha,m)=&\left\{u\in W^{2,p}_{loc}(\R^{N+1}_+):\   y^{\alpha} D^2u,\ y^\frac{\alpha}{2} \nabla u,\ u\in L^p_m\right\}.
\end{align*}
\medskip 

We start by considering, preliminarily, the case $d=0$ and we add, accordingly,  a Neumann  boundary condition 
$$\ds \lim_{y\to 0}y^{\frac c\gamma}D_yu=0$$
 in  integral  form (see  Proposition \ref{Trace D_yu in W}) by defining
\begin{align*}
	W^{2,p}_{\mathcal N}(\alpha,\alpha,m)&=\{u \in W^{2,p}(\alpha,\alpha,m):\  y^{\alpha-1}D_yu\ \in L^p_m\}.
\end{align*}

The following result characterizes the generation properties, the maximal regularity and the   domain of $\mathcal L$ in the special case $d=0$. It has been proved by the author in \cite{MNS-Degenerate-Half-Space} by constructing, using tools from vector-valued harmonic analysis and  Fourier multipliers, a resolvent  $(\lambda-\mathcal L)^{-1}$ of $\mathcal L$ for $\lambda$ in a suitable sector $\omega+\Sigma_{\phi}$, with $\omega \in \R$ and by proving that the family  $\lambda (\lambda-\mathcal L)^{-1}$ is $\mathcal  R$-bounded on $\mathcal B\left(L^p_m\right)$.  We refer the reader to \cite{MNS-Degenerate-Half-Space} for the proof and any further details.

\begin{teo} \label{Teo alpha Neum}{\em \cite[Theorem 5.6,  Corollary 5.7, Theorem 6.3]{MNS-Degenerate-Half-Space}}
	Let $\alpha\in\R$ such that $\alpha<2$ and
	$$\alpha^{-} <\frac{m+1}p<\frac{c}{\gamma}+1-\alpha.$$   
	Then the operator
	\begin{align*}
		\mathcal L&=y^{\alpha}\mbox{Tr }\left(AD^2\right)+y^{\alpha-1} c  D_y
	\end{align*}
endowed with domain   $W^{2,p}_{\mathcal N}\left(\alpha,\alpha,m\right)$  generates a bounded analytic semigroup  in $L^p_m$ which has maximal regularity. Moreover the set  $C_c^\infty (\R^{N})\otimes \mathcal D$ defined in \eqref{defDcore} is a core for $\mathcal L$. 
Finally, the estimate
\begin{equation*} 
\|y^\alpha D_{x_i x_j} u\|_{L^p_{m}} +\|y^\alpha D_{yy} u\|_{L^p_{m}}+\|y^\alpha D_{x_iy} u\|_{L^p_{m}}+\|y^{\alpha-1} D_{y} u\|_{L^p_{m}}\leq C\| \mathcal Lu\|_{L^p_{m}}
\end{equation*}
holds for every $u \in W^{2,p}_{\mathcal{N}}(\alpha,\alpha,m)$
\end{teo}
\medskip 
When $c\neq 0$ and $d\in\R$, we can   impose an oblique derivative boundary condition 
$$\lim_{y\to 0}y^{\frac c\gamma}\, v \cdot \nabla u=0$$
in integral form (see Proposition \ref{Equiv Wv alpha eq})  by defining 
\begin{align*}
	W^{2,p}_{v}(\alpha, \alpha,m)\stackrel{def}{=}\{u \in W^{2,p}(\alpha, \alpha,m):\  y^{\alpha-1}v\cdot \nabla u \in L^p_m\}.
\end{align*}
Recalling Definition \ref{convention Wv Wn alpha}, to shorten the notation we also write $W^{2,p}_{(0,0)}(\alpha, \alpha,m)=W^{2,p}_{\mathcal N}(\alpha, \alpha,m)$.
\medskip

We transform $\mathcal L$  into a similar operator with $d=0$ and Neumann boundary condition. 
Indeed,  we  use the map $S_{0,\omega}$ of Section \ref{Section Degenerate}   defined in \eqref{Tran def} with  $\beta=0$ and  $\omega=-\frac d {c}$, namely
\begin{align*}
	S_{0,-\frac d {c}}\, u(x,y)&:=u\left(x-\frac d {c} y,y\right),\quad (x,y)\in\R_+^{N+1}.
\end{align*}
We recall that, by Proposition \ref{Isometry shift} and 
Corollary \ref{Isometry action 2 spe},  $S_{0,-\frac d {c}}$ is an isometry of  $L^p_{m}$  and for every
$u\in W^{2,1}_{loc}\left(\R^{N+1}_+\right)$ one has
	\begin{align*}
	{S_{0,-\frac d {c}}}^{-1} &\Big(y^{\alpha}\mbox{Tr }\left(AD^2\right)+y^{\alpha-1}\left(v,\nabla\right)\Big){S_{0,-\frac d {c}}}u=y^{\alpha}\mbox{Tr }\left(\tilde AD^2u\right)+y^{\alpha-1}c D_yu
\end{align*}
where 
\begin{align*}
	\tilde A=	\left(
	\begin{array}{c|c}
		\tilde Q  & { \tilde q}^t \\[1ex] \hline
		\tilde q& \gamma
	\end{array}\right),\qquad \tilde Q=Q-\frac{2}c&q\otimes  d +\frac{\gamma}{c^2}d\otimes d,\qquad \tilde q=q-\frac{\gamma}c d.
\end{align*}

We can then deduce the following result.

\begin{teo}\label{Teo alpha Oblique} 
	Let $v=(d,c)\in\R^{N+1}$ with $d=0$ if $c=0$, and let $\alpha\in\R$ such that $\alpha<2$ and
	$$\alpha^{-} <\frac{m+1}p<\frac{c}{\gamma}+1-\alpha.$$   
	Then the operator
	\begin{align*}
		\mathcal L&=y^{\alpha}\mbox{Tr }\left(AD^2\right)+y^{\alpha-1}\left(v,\nabla\right)
	\end{align*}
   generates a bounded analytic semigroup  in $L^p_m$ which has maximal regularity. Moreover 
   \begin{align*}
   	D(\mathcal L)=W^{2,p}_{v}\left(\alpha,\alpha,m\right)
   \end{align*} 
and the set  $ \mathcal C_v$ defined in \eqref{defC oblique} is a core for $\mathcal L$.
\end{teo}
{\sc{Proof.}}
 The claim for $c=0$ or $d=0$ is just  Theorem \ref{Teo alpha Neum}. Let us suppose $c\neq 0$. According to the discussion above the isometry $S_{0,-\frac d {c}}$ of $L^p_m$
transforms  $\mathcal L$ into
\begin{align*}
	\tilde {\mathcal L}=y^{\alpha}\mbox{Tr }\left(\tilde AD^2\right)+y^{\alpha-1}c D_y,
\end{align*} 
The statement on generation and maximal regularity is therefore a translation to  $\mathcal L$ and in $L^p_m$ of the results of Theorem \ref{Teo alpha Neum} for $\tilde{\mathcal L}$.
Also, using Proposition \ref{Equi obli Neumm sob alpha},   one has 
\begin{align*}
	D(\mathcal L)=S_{0,-\frac d {c}}\left(D\left(\tilde {\mathcal L}\right)\right)=S_{0,-\frac d {c}}\left(W^{2,p}_{\mathcal N}\left(\alpha,\alpha,m\right)\right)=W^{2,p}_{v}\left(\alpha,\alpha,m\right).
\end{align*}
\qed

\begin{cor}\label{Oblique cor1}
Under the assumptions of the previous theorem, the estimate
	\begin{align}\label{elliptic regularity oblique} 
		\|y^\alpha D_{x_i x_j} u\|_{L^p_{m}} +\|y^\alpha D_{yy} u\|_{L^p_{m}}+\|y^\alpha D_{x_iy} u\|_{L^p_{m}}+\|y^{\alpha-1} v\cdot \nabla u\|_{L^p_{m}}\leq C\| \mathcal Lu\|_{L^p_{m}}
	\end{align}
\mbox{}\\[-1.5ex]	holds for every $u \in W^{2,p}_{v}(\alpha,\alpha,m)$ (if $c=0$ replace $y^{\alpha-1} v\cdot \nabla u$ with $y^{\alpha-1} D_yu$).
\end{cor}
{\sc{Proof.}}
By Theorem \ref{Teo alpha Oblique}  the above inequality holds if $\|  u\|_{L^p_{m}}$ is added to the right hand side.
Applying it  to $u_\lambda (x,y)=u(\lambda x, \lambda y)$, $\lambda >0$ we obtain
\begin{align*}
	\|y^\alpha  D_{x_i x_j} u\|_{L^p_{m}}+\|y^\alpha  D_{x_i y} u\|_{L^p_{m}} +\|y^\alpha D_{yy} u\|_{L^p_{m}}+\|y^{\alpha-1}v\cdot\nabla  u\|_{L^p_{m}}\leq C\left(\| \mathcal L u\|_{L^p_{m}}+\lambda^{\alpha -2}\| u\|_{L^p_{m}}\right)
\end{align*}
and the proof follows letting $\lambda \to \infty$.\qed 

The following corollary  enlightens the role of the  Neumann  and of the oblique boundary conditions. 
\begin{cor}\label{Oblique cor2}
	Under the assumptions of the previous theorem one has 
\begin{align*} 
	D(\mathcal L)=&\left\{u \in  W^{2,p}_{loc}(\R^{N+1}_+): u,\ y^{\alpha}\Delta_xu,\  y^{\alpha}\gamma D_{yy}u+y^{\alpha-1}v\cdot \nabla u \in L^p_m\text{\;\;and\;\;}\lim_{y\to 0}y^{\frac c\gamma} v\cdot \nabla u =0\right\}
\end{align*}
(when $v=0$, replace \;$\ds\lim_{y\to 0}y^{\frac c\gamma} v\cdot \nabla u=0$\; with  \;$\ds\lim_{y\to 0} D_y u=0$).
\end{cor}
{\sc{Proof.}} The proof follows from Theorem \ref{Teo alpha Oblique}  and Propositions \ref{Trace D_yu in W} and \ref{Equiv Wv alpha eq}.\qed

 \section{The operator with Dirichlet boundary conditions}\label{Section DIrichlet alpha}

 In this section we add a potential term to the operator of Section \ref{Sect DOm obliqu}  and study  the  operator  
\begin{align*}
	\mathcal L&=y^{\alpha}\mbox{Tr }\left(AD^2\right)+y^{\alpha-1}\left(v,\nabla\right)-by^{\alpha-2}
\end{align*}
 with $v=(d,c)\in\R^{N+1}$ and $b\in\R$.  We endow $\mathcal L$ under Dirichlet boundary conditions, in the sense specified below.
 
We always assume the hypotheses in  Assumptions   \ref{assumption} and \ref{Assumption D}; we also recall that, as in  Definition \ref{Definition operator}, $\mathcal L$ can be written equivalently as
 \begin{align*}
 	\mathcal L&=y^{\alpha}\mbox{Tr }\left(QD^2_x\right)+2y^{\alpha}\left(q,\nabla_xD_y\right)+y^{\alpha}\gamma  L_y+y^{\alpha-1}\left(d,\nabla_x\right)
 \end{align*}
 where $A=	\left(
 \begin{array}{c|c}
 	Q  & { q}^t \\[1ex] \hline
 	q& \gamma
 \end{array}\right)$ and $$L_y=D_{yy}+\frac{c/\gamma}{y}D_y-\frac{b/\gamma}{y^2}$$ is the operator defined in \eqref{def 1-d op} with parameters $\frac b\gamma,\ \frac c\gamma$.

 The equation $L_yu=0$ has solutions $y^{-s_1}$, $y^{-s_2}$ where $s_1,s_2$ are the roots of the indicial equation $f(s)=-s^2+(c/\gamma-1)s+b/\gamma=0$ given by
 
 \begin{equation} \label{defs gamma}
 	s_1:=\frac{\frac c \gamma-1}{2}-\sqrt{D},
 	\quad
 	s_2:=\frac{ \frac c \gamma-1}{2}+\sqrt{D}
 \end{equation}
 where
 \bigskip
 \begin{equation*} 
 	D:=
 	\frac b\gamma+\left(\frac{\frac c \gamma -1}{2}\right)^2
 \end{equation*}
 is supposed to be nonnegative. 
 When $b=0$, then $\sqrt D=|c/\gamma-1|/2$ and $s_1=0, s_2=c/\gamma-1$ for $c/\gamma \ge 1$ and $s_1=c/\gamma-1, s_2=0$ for $c/\gamma<1$.
 
\begin{os} 
 	All the results of this section will be valid, with obvious changes,  also  in  $\R_+$ for the $1$d operators $y^{\alpha} L_y$ changing (when it appears in the various conditions on the parameters) $\alpha^{-}$ to $0$ (see also Remark  \ref{Os Sob 1-d}).  We  also refer to   \cite{met-calv-negro-spina, MNS-Sharp, MNS-Grad, MNS-Max-Reg,  Negro-Spina-Asympt} for the analogous results concerning the $Nd$ version of $L_y$.
 \end{os}

 \medskip
 A  multiplication operator transforms $\mathcal L$  into an operator of the form $y^{\alpha}\mbox{Tr }\left(AD^2\right)+y^{\alpha-1}\left(v,\nabla\right)$  and allows  to transfer the results of Section \ref{Sect DOm obliqu} to this  situation. 
 Indeed,  we  use the map defined in  \eqref{Gen Kelvin def} of Section \ref{Section Degenerate}
 \begin{align}\label{Gen Kelvin def beta 0}
 	T_{k,0\,}u(x,y)&:=y^ku(x,y),\quad (x,y)\in\R^{N+1}_+
 \end{align}
 for a suitable choice of $k$ and with $\beta=0$.
 We recall that, by Propositions \ref{Isometry action der}, \ref{Isometry action} and Corollary \ref{Isometry action2},  $T_{k,0\,}$ maps isometrically  $L^p_{\tilde m}$ onto $L^p_m$  where 
 $ \tilde m=m+kp$ and for every
 $u\in W^{2,1}_{loc}\left(\R^{N+1}_+\right)$ one has
\begin{align*}
	& T_{k,0\,}^{-1} \Big(\mathcal L\Big)T_{k,0\,}u=\tilde {\mathcal L}u
\end{align*}
where 
\begin{align}
	\label{tilde b section}
\nonumber \tilde{\mathcal L}=y^{\alpha}\mbox{Tr }\left(AD^2\right)+y^{\alpha-1}\left(\tilde d,\nabla_x\right)+y^{\alpha-1}\tilde cD_y-\tilde by^{\alpha-2},\\[1.5ex]
\tilde d=2kq+d,\qquad \tilde c=c+2k\gamma,\qquad \tilde b=b-k\left(c+(k-1)\gamma\right).
\end{align}
Equivalently  we can  write
\begin{align*}
	\tilde{\mathcal L}	&=y^{\alpha}\mbox{Tr }\left(QD^2_xu\right)+2y^{\alpha}\left(q,\nabla_xD_yu\right)+y^{\alpha}\gamma  \tilde L_yu+y^{\alpha-1}\left(2kq+d,\nabla_xu\right)
\end{align*}
 where $\tilde { L}_y=D_{yy}+\frac{\tilde c/\gamma}{y}D_y-\frac{\tilde b/\gamma}{y^2}$. 
 Moreover the discriminant $\tilde D$ and the parameters $\tilde s_{1,2}$ of $\tilde L_y$ are given by
 \begin{align}\label{tilde D gamma section}
 	\tilde D&=D, \quad\tilde s_{1,2}=s_{1,2}+k.
 \end{align}
 
  Choosing $k=-s_1$ and recalling the definition of $s_1$,  we get 
  \begin{align}\label{tilde c gamma}
  	\tilde b=0,  \qquad \tilde c=\gamma\left(\frac c\gamma-2s_1\right)=\gamma\left(1+2\sqrt D\right)
  \end{align} 
 and therefore  
 \begin{align}\label{tilde l gamma}
 	T_{-s_1,0}^{-1}\Big(y^{\alpha}\mbox{Tr }\left(AD^2\right)+y^{\alpha-1}\left(v,\nabla\right)-by^{\alpha-2}\Big)T_{-s_1,0}=y^{\alpha}\mbox{Tr }\left(AD^2\right)+y^{\alpha-1}\left(w,\nabla\right)
 \end{align}
 where $w:=\left(\tilde d, \tilde{c}\right)=\left(d-2s_1q, c-2s_1\gamma\right)$. Moreover from \eqref{tilde c gamma} one has  \begin{align*}
 	w=v-2s_1\left(q,\gamma\right)=\left(d-2s_1q, \gamma(1+2\sqrt D)\right).
 \end{align*}

We can now  derive the following result.
 \begin{teo} \label{complete}	Let $\alpha\in\R$ such that $\alpha<2$ and
 	$$s_1+ \alpha^-<\frac{m+1}p<s_2+2-\alpha.$$   
 	Then, under Assumptions   \ref{assumption} and \ref{Assumption D}, the operator
 	\begin{align*}
 		\mathcal L&=y^{\alpha}\mbox{Tr }\left(AD^2\right)+y^{\alpha-1}\left(v,\nabla\right)-by^{\alpha-2}
 	\end{align*} 
  generates a bounded analytic semigroup  in $L^p_m$ which has maximal regularity. Moreover,
 	\begin{equation}
 		\label{dominioTrasf alpha}
 		D(\mathcal L)
 		=y^{-s_1}W^{2,p}_{w}\left(\alpha,\alpha,m-s_1p\right),\qquad w=v-2s_1\left(q,\gamma\right).
 	\end{equation}	
 	
 
 \end{teo}
 {\sc Proof.}  According to the discussion above the map
 $T_{-s_1,0}:L^p_{m-s_1p}\to L^p_m$
 transforms  $\mathcal L$ into
 \begin{align}\label{eq complete1}
 	 \tilde {\mathcal L}=y^{\alpha}\mbox{Tr }\left(AD^2\right)+y^{\alpha-1}\left(w,\nabla\right),\qquad w=\left(-2s_1q+d, c-2s_1\gamma\right):=\left(\tilde d,\tilde c\right).
 \end{align} 

We observe now that, recalling \eqref{defs gamma}, the hypothesis  $s_1+ \alpha^-<\frac{m+1}p<s_2+2-\alpha$ is equivalent to $ \alpha^-<\frac{m-ps_1+1}p<\frac{\tilde c}{\gamma}+1-\alpha$. Moreover  Assumption  \ref{Assumption D} implies, recalling \eqref{tilde c gamma}, that  $\tilde c=\gamma(1+2\sqrt D)\geq \gamma>0$. 
 The statement on generation and maximal regularity is therefore a translation to  $\mathcal L$ and in $L^p_m$ of the results of Theorem \ref{Teo alpha Oblique} for  $y^{\alpha}\mbox{Tr }\left(AD^2\right)+y^{\alpha-1}\left(w,\nabla\right)$ in $L^p_{m-s_1p}$. 
 
 Also  one has 
 \begin{align*}
 	D(\mathcal L)=T_{-s_1,0}\left(D\left(\tilde {\mathcal L}\right)\right)=y^{-s_1}W^{2,p}_{w}\left(\alpha,\alpha,m-s_1p\right).
 \end{align*}
 \qed

  \begin{os}\label{oss ellipric regularity}\mbox{}
  	\begin{itemize}
 \item[(i)] Recalling Definition \eqref{Definition W_v alpha} and by a straightforward calculation, equality \eqref{dominioTrasf alpha} says   that $u \in y^{-s_1}W^{2,p}_{w}\left(\alpha,\alpha,m-s_1p\right)$ if and only if   all functions
\begin{align*}
 	u,\hspace{1.3ex}  y^{\alpha}D_{x_i x_j}u,\hspace{1.3ex} y^\alpha\left(\nabla_xD_yu+s_1\frac{\nabla_xu}y\right),\hspace{1.3ex} y^{\frac{\alpha}{2}} \nabla_xu,\hspace{1.3ex} y^{\frac \alpha 2}\left (D_y u+s_1\frac{u}{y} \right ),\\[1ex]
y^\alpha\left(D_{yy}u+2s_1\frac{D_yu}{y}-\left(s_1-s_1^2\right)\frac{u}{y^2}\right),
\hspace{2ex} y^{\alpha}\left( \frac{w\cdot \nabla u}y+s_1\tilde c\frac{u}{y^2}\right)
 \end{align*} 
  belong to $L^p_m$ ( $\tilde c=w\cdot e_{N+1}=c-2s_1\gamma$).  However, in the range of parameters of Theorem \ref{complete}, one cannot deduce, in general, that $y^{\alpha-1}w\cdot\nabla u$ and $y^{\alpha}D_{yy}u$ belong to $L^p_m$, as one can check on functions like $y^{-s_1}u\left(x-\frac{\tilde d}{\tilde c}y\right)$, $u\in C_c^\infty(\R^N)$, near $y=0$.  This is however possible in the special case of Corollary \ref{Cor-Rellich}.
\item[(ii)] 	Unlike Theorem \ref{Teo alpha Oblique}, the above remark shows that  one cannot in general estimate any singular terms which compose $\mathcal L$ as in \eqref{elliptic regularity oblique}. 
 	Nevertheless, the estimate involving the $x$-derivatives
 	\begin{equation*} 
 		\|y^{\alpha}D_{x_ix_j} u\|_{L^p_m} \leq C \|\mathcal L u\|_{L^p_m}
 	\end{equation*}
 and, by difference, also \; $\|y^{\alpha}\mathcal L u-y^\alpha \mbox{Tr }(QD^2_x u) \|_{L^p_m} \leq C \|\mathcal L u\|_{L^p_m}$,\;
 always	hold for every $u \in D(\mathcal L)$. This follows 
since, by Corollary \ref{Oblique cor1},  the similar statement holds for $\tilde {\mathcal L}$ in $L^p_{m-s_1p}$ and,  by Propositions \ref{Isometry action der} and \ref{Isometry action}, 
$$T_{-s_1,0}^{-1}\left( y^{\alpha} D_{x_ix_j}\right) T_{-s_1,0}=y^{\alpha}D_{x_ix_j}.
$$
  \end{itemize} \end{os}

\medskip 
 The following corollary explains why we use the term Dirichlet boundary conditions.

 \begin{cor} \label{cor1}
 	Let $\alpha\in\R$ such that $\alpha<2$ and   $s_1+ \alpha^-<\frac{m+1}p<s_2+2-\alpha$ and let $\mathcal A$ be the operator
 	\begin{align*}
 		\mathcal A:&=\gamma y^{\alpha}D_{yy}+y^{\alpha-1} v_q \cdot \nabla-by^{\alpha-2},\qquad v_q:=\left(d-2s_1q,c\right).
 	\end{align*}
 	\begin{itemize}
 		\item[(i)] If $D>0$ then 
 		\begin{align*} 
 	D(\mathcal L)=&\left\{u \in  W^{2,p}_{loc}(\R^{N+1}_+): u,\ y^{\alpha}\Delta_xu,\ \mathcal Au \in L^p_m\text{\;\;and\;\;}\lim_{y\to 0}y^{s_2} u=0\right\}
 	\end{align*}
 		\item[(ii)] If $D=0$ then $s_1=s_2$ and
 	\begin{align*} 
 		D(\mathcal L)=&\left\{u \in  W^{2,p}_{loc}(\R^{N+1}_+): u,\ y^{\alpha}\Delta_xu,\ \mathcal Au\in L^p_m\text{\;\;and\;\;}\lim_{y\to 0}y^{s_2} u\in\C\right\}
 	\end{align*}
 	\end{itemize} 
 	In both cases the graph   norm and $$\|u\|_{L^p_m}+\|y^{\alpha}\Delta_x u\|_{L^p_m}+\|\mathcal Au\|_{L^p_m}$$ are equivalent on $D(\mathcal L)$. 
 \end{cor}
 {\sc Proof.} Let us prove claim (i). By Theorem \ref{complete} and  Remark \ref{oss defi W_v} (iv) we have 
 \begin{equation*}
 	D(\mathcal L)
 	=T_{-s_1,0}\left(W^{2,p}_{w}\left(\alpha,\alpha,m-s_1p\right)\right)=T_{-s_1,0}\left(W^{2,p}_{\frac w\gamma}\left(\alpha,\alpha,m-s_1p\right)\right)
 \end{equation*}	
 where
 \begin{equation*}
 	w=(d,c)-2s_1\left(q,\gamma\right):=\left(\tilde d,\tilde c\right).
 \end{equation*}	
We apply  Proposition \ref{Equiv Wv alpha dir}  with  $\frac{\tilde{c}}{\gamma}=\frac{c}{\gamma}-2s_1=1+2\sqrt{D}> 1$ in place of $c$ thus obtaining
 \begin{align*} 
 	W^{2,p}_{w}(\alpha, \alpha, m-s_1p)=&\left\{v \in  W^{2,p}_{loc}(\R^{N+1}_+): v,\ y^{\alpha}\Delta_xv,\   	\tilde{\mathcal A}v \in L^p_{m-s_1p}\text{\;\;and\;\;}\lim_{y\to 0}y^{\frac{\tilde c}{\gamma}-1} v=0\right\}
 \end{align*}
 where
 \begin{align*}
 	\tilde{\mathcal A}:&=y^{\alpha}D_{yy}+y^{\alpha-1} \frac \omega\gamma  \cdot \nabla.
 \end{align*}
 The required claim then follows from the previous equalities after noticing that, by Proposition \ref{Isometry action}, 
 \begin{align*}
 	&T_{-s_1,0}^{-1}\left( y^{\alpha} \Delta_x\right) T_{-s_1,0}=y^{\alpha}\Delta_x,\qquad  T_{-s_1,0}^{-1}\left(\mathcal A \right)T_{-s_1,0}=\gamma \tilde{\mathcal A}.
 \end{align*} 
 and that  for any $u\in y^{-s_1}W^{2,p}_{w}\left(\alpha,\alpha,m-s_1p\right)$, setting $v=y^{s_1}u$, one has, recalling \eqref{defs},
 \begin{align*}
 	y^{\frac{\tilde c}{\gamma} -1}v=y^{\frac{\tilde c}{\gamma}-1+s_1}u=y^{2\sqrt D+s_1}u=y^{s_2 }u.
 \end{align*} 
 Claim (ii) follows similarly.
 \qed

 In the following Corollary, in continuity with Remark \ref{oss ellipric regularity} (ii), we shows that in certain range of parameters, we can improve the elliptic regularity of the operator $\mathcal L$ by estimating, in addition to $y^{\alpha}D_{x_ix_j}u$, other different terms which compose $\mathcal L$.  Specifically if $\frac{m+1}p>s_1+1-\alpha$
 we can estimate  
 \begin{align*}
 	y^{\alpha} \gamma L_yu,\quad y^\alpha D_{x_iy} u,\quad y^{\alpha-1}  \nabla_x u
 \end{align*}
 whereas in the smaller  range  $\frac{m+1}p>s_1+2-\alpha$, we reach the best elliptic regularity where, recalling the definition of $W^{2,p}_{\mathcal R}(\alpha,\alpha,m)$ in Section \ref{section W_r}, $D(\mathcal L)$ consists of the functions for which any single terms of $\mathcal L$ is in $L^p_m$.
 
 \begin{cor} \label{Cor-Rellich}	Let $\alpha\in\R$ such that $\alpha<2$. 
 	\begin{itemize}
 	\item[(i)] If both the condition  $s_1+ \alpha^-<\frac{m+1}p<s_2+2-\alpha$ and  $\frac{m+1}p>s_1+1-\alpha$ hold then $D(\mathcal L)
 	=y^{-s_1}W^{2,p}_{\mathcal N}\left(\alpha,\alpha,m-s_1p\right)$ and 
 	\begin{equation*} 
 		\|y^\alpha D_{x_i x_j} u\|_{L^p_{m}} +\| y^{\alpha} \gamma L_yu\|_{L^p_{m}}+\|y^\alpha D_{x_iy} u\|_{L^p_{m}}+\|y^{\alpha-1}  \nabla_x u\|_{L^p_{m}}\leq C\| \mathcal Lu\|_{L^p_{m}}
 	\end{equation*}
 	where $ y^{\alpha} \gamma L_y=y^\alpha \gamma D_{yy} +cy^{\alpha-1}D_y-by^{\alpha-2}$.
 	\item[(ii)] If $s_1+2-\alpha<\frac{m+1}p<s_2+2-\alpha$ then $D(\mathcal L)=W^{2,p}_{\mathcal R}(\alpha,\alpha,m)$.
 	\end{itemize}
 	
 \end{cor}
 {\sc Proof.} Following the notation of the proof of Theorem  \ref{complete}, we consider the operator $\tilde {\mathcal L}$  of \eqref{eq complete1} on $L^p_{m-s_1p}$.
  Claim (i) then follows since by Proposition \ref{Hardy Rellich Sob}, 
$W^{2,p}_{w}\left(\alpha,\alpha,m-s_1p\right)=W^{2,p}_{\mathcal N}\left(\alpha,\alpha,m-s_1p\right)$. Moreover if $u\in D(\mathcal L)=y^{-s_1}W^{2,p}_{\mathcal N}\left(\alpha,\alpha,m-s_1p\right)$, then $v=y^{s_1}u$ satisfies  by Proposition \ref{Hardy Rellich Sob} and by \eqref{elliptic regularity oblique}  of Theorem \ref{Teo alpha Oblique}
 \begin{equation*} 
 \|y^{\alpha-1} \nabla_x v\|_{L^p_{m-s_1p}}\leq C\| \tilde{\mathcal L}v\|_{L^p_{m-s_1p}}
 \end{equation*}
which under the isometry $T_{-s_1,0}:L^p_{m-s_1p}\to L^p_m$ translates, recalling that $T_{-s_1,0}^{-1}\left( y^{\alpha} \nabla_x\right) T_{-s_1,0}=y^{\alpha}\nabla_x$ and $T_{-s_1,0}^{-1}\left( \mathcal L\right) T_{-s_1,0}=\tilde{\mathcal L}$, into
 \begin{equation*} 
	\|y^{\alpha-1} \nabla_x u\|_{L^p_{m}}\leq C\|{\mathcal L}u\|_{L^p_{m}}.
\end{equation*}
Then by difference, using the equivalent norm of Corollary \ref{cor1}, the last estimate proves that the required inequality holds if $\|  u\|_{L^p_{m}}$ is added to the right hand side. This term can although be eliminated by performing the same scaling argument used in the proof of Theorem \ref{Teo alpha Oblique}.

 To prove  (ii) we observe preliminarily that  $s_1+ 2-\alpha>s_1+\alpha^- $, since  $\alpha<2$. By  Theorem \ref{complete} and Proposition \ref{RN}
 $$D(\mathcal L)=y^{-s_1}\left( W^{2,p}_{w}\left(\alpha,\alpha,m-s_1p\right)\right)=W^{2,p}_{\mathcal R}\left(\alpha,\alpha,m\right)$$
 under the assumption $\frac{m-ps_1+1}p>2-\alpha$ which is equivalent to $s_1+ 2-\alpha<\frac{m+1}p$.
 \qed

 \medskip 
 
 Observe that the condition $\frac{m+1}p>s_1+1-\alpha$ in the previous corollary is necessary for the integrability  of the mixed derivatives of functions like $y^{-s_1}u(x)$, $u\in C_c^\infty(\R^N)$, near $y=0$.
 
 
 \medskip\medskip 
 
 The above results apply, when $b=0$, also to the operator
 \begin{align*}
 	\mathcal L&=y^{\alpha}\mbox{Tr }\left(AD^2\right)+y^{\alpha-1}\left(v,\nabla\right).
 \end{align*}

 \begin{os}\label{oss uniquenss DirNeu}
 	
When $c>\gamma $, so that $s_1=0$, $s_2=\frac c\gamma -1 > 0$, by \eqref{dominioTrasf alpha} the operator $\mathcal L$  coincides with the one of Theorems \ref{Teo alpha Neum} and \ref{Teo alpha Oblique} since $D(\mathcal L)=W^{2,p}_{v}\left(\alpha,\alpha,m\right)$. Moreover by Proposition \ref{Equiv Wv alpha diff dir} and Corollary \ref{cor1} one has
\begin{align*}
	D(\mathcal L)	&=\left\{u \in  W^{2,p}_{loc}(\R^{N+1}_+): u,\ y^{\alpha}\Delta_xu,\ y^{\alpha}\gamma D_{yy}u+y^{\alpha-1}v\cdot \nabla u \in L^p_m\text{\;\;and\;\;}\lim_{y\to 0}y^{\frac c\gamma -1} u=0\right\}.
\end{align*}
Therefore, recalling Corollary \ref{Oblique cor2} the oblique derivative and the Dirichlet boundary conditions are, in the range $c>\gamma $,  equivalent since they lead to the same operator.  
 \end{os}
 
 On the other hand, when $c<\gamma $, so that $s_1=\frac c\gamma -1 \neq 0$, $s_2=0$, we can construct a realization of $\mathcal L$  different from that of Theorems \ref{Teo alpha Neum} and \ref{Teo alpha Oblique}. 
 
 \begin{cor}\label{CorBdalpha}
 	Let   $c<\gamma$  and $\frac c\gamma -1+ \alpha-<\frac{m+1}p<2-\alpha.$ Then $\mathcal L=y^{\alpha}\mbox{Tr }\left(AD^2\right)+y^{\alpha-1}\left(v,\nabla\right)$ with domain
 	 		\begin{align*} 
 		D(\mathcal L)=&\left\{u \in  W^{2,p}_{loc}(\R^{N+1}_+): u,\ y^{\alpha}\Delta_xu,\ y^{\alpha}\gamma D_{yy}u+y^{\alpha-1}v_q\cdot \nabla u \in L^p_m\text{\;\;and\;\;}\lim_{y\to 0} u=0\right\},
 	\end{align*}
 	where $v_q:=\left(d-2s_1q,c\right)\in\R^{N+1}$,
 	generates a bounded analytic semigroup  in $L^p_m$ which has maximal regularity.
 \end{cor}
 {\sc Proof. } This follows from Corollary \ref{cor1} (i), since $s_1=\frac c\gamma -1$ and $s_2=0$.
 \qed
 Note that the generation interval  $\frac c\gamma-1+ \alpha^-<\frac{m+1}p<2-\alpha$ under Dirichlet boundary conditions,   is larger than $ \alpha^-<\frac{m+1}p<\frac c\gamma+1-\alpha$  given by  Theorems \ref{Teo alpha Neum} and \ref{Teo alpha Oblique} for Neumann and oblique  boundary conditions.
 

 	\begin{os}
 	Let us explain what happens in Theorem \ref{complete} if we choose $-k$ in \eqref{tilde b section} as the second root $s_2$ instead of $s_1$. Proceeding similarly, one proves an identical result under the condition
 \begin{equation} \label{nonunique}
 	s_2+\alpha^-< \frac{m+1}{p} <s_1+2-\alpha.
 \end{equation}
 However this requires the assumption $s_2<s_1+2-\alpha$ which is not always satisfied. When \eqref{nonunique} holds this procedure leads to a different operator. For further details about  different realizations of $\mathcal L$  and about  uniqueness questions  we refer also to \cite[Section 9.2]{MNS-CompleteDegenerate}.
 	\end{os}


\section{Consequences for more general operators}\label{Consequenses}
%
In this section we  deduce generation and domain properties in  $L^p_m$ for the more general operators
	\begin{align}\label{operator alpha diff def}
		\mathcal L=y^{\alpha_1}\mbox{Tr }\left(QD^2_x\right)+2y^{\frac{\alpha_1+\alpha_2}{2}}q\cdot \nabla_xD_y+\gamma y^{\alpha_2}D_{yy}+y^{\frac{\alpha_1+\alpha_2}{2}-1}\left(d,\nabla_x\right)+c y^{\alpha_2-1}D_y-by^{\alpha_2-2}
	\end{align}
where possibly different powers $y^{\alpha_1}$, $y^{\alpha_2}$ appear in front respectively  of the $x$ and $y$ derivatives. Here $\alpha_1,\alpha_2\in\R$ such that   $\alpha_2<2$, $\alpha_2-\alpha_1<2$. We keep the Assumptions  \ref{assumption} and \ref{Assumption D} and we recall that $\mathcal L$ can be written equivalently in the compact form 
	\begin{align*}
		\mathcal L&=y^{\alpha_1}\mbox{Tr }\left(QD^2_x\right)+2y^{\frac{\alpha_1+\alpha_2}2}\left(q,\nabla_xD_y\right)+\gamma y^{\alpha_2} L_y+y^{\frac{\alpha_1+\alpha_2}{2}-1}\left(d,\nabla_x\right)
	\end{align*}
	where $$ L_y=\gamma D_{yy}+\frac {c/\gamma}y D_y-\frac{b/\gamma}{y^{2}}$$ is the operator defined in \eqref{def 1-d op}. 
	
	 In contrast to the case $\alpha_1=\alpha_2$, some complications arise due to the different weights $y^{\alpha_1}$, $y^{\alpha_2}$ which appear in the $x$ and $y$ directions. This reflects into  a distortion correction  which depends on  the coefficient
\begin{align*}
	\beta_{\alpha}\overset{def}=\frac{\alpha_1-\alpha_2}{2}
\end{align*} 
which appears into the characterization of the domain of $\mathcal L$ (see the definition of $W^{2,p}_{v}(\alpha_1,\alpha_2,m)$  in Section \ref{sobolev oblique alpha diff section}). This complication  does not appears, obviously,  when $\alpha_1=\alpha_2$, where $\beta_\alpha=0$, but  also when $d=0$ (see Remark \ref{oss defi W_v}). \\

The isometry of Section \ref{Section Degenerate} transforms $\mathcal L$  into a similar operator with $\alpha_1=\alpha_2$  and allows  to transfer the results of Sections \ref{Sect DOm obliqu} and \ref{Section DIrichlet alpha} to this  situation.

Indeed,  we  use the map  defined in  \eqref{Gen Kelvin def}  with $k=0$ and $\beta={\beta_\alpha}$ namely
\begin{align}\label{Gen Kelvin def k 0}
T_{0,{\beta_\alpha}}u(x,y)&:=|{\beta_\alpha}+1|^{\frac 1 p}u(x,y^{{\beta_\alpha}+1}),\quad (x,y)\in\R^{N+1}_+,\qquad  {\beta_\alpha}=\frac{\alpha_1-\alpha_2}{2}.
\end{align}
We recall that, by Propositions \ref{Isometry action der},  \ref{Isometry action} and Corollary \ref{Isometry action2},  $T_{0,{\beta_\alpha}}$ maps isometrically  $L^p_{\tilde m}$ onto $L^p_m$  where 
$ \tilde m=\frac{m-{\beta_\alpha}}{{\beta_\alpha}+1}$ and
  transforms  $\mathcal L$ into   
\begin{align*}
		\tilde{\mathcal L}:=	T_{0,{\beta_\alpha}}^{-1} \,\mathcal L\, T_{0,{\beta_\alpha}}&=y^{\alpha}\mbox{Tr }\left(\tilde AD^2\right)+y^{\alpha-1}\left(\tilde v,\nabla\right)- by^{\alpha-2}
\end{align*}
where
\begin{align}\label{coeff tilde L eq}
\nonumber 	\alpha=\frac{\alpha_1}{{\beta_\alpha}+1},\qquad &\tilde A:=	\left(
\begin{array}{c|c}
	Q  & {\tilde q}^t \\[1ex] \hline
	\tilde q& \tilde \gamma
\end{array}\right),\quad 
\tilde q=({\beta_\alpha}+1)q,\qquad \tilde \gamma=({\beta_\alpha}+1)^2\gamma,\\[1.5ex]
&\tilde v:=( d, \tilde c),\qquad \tilde c=(c+{\beta_\alpha}\gamma)({\beta_\alpha}+1).
\end{align}
 As in  Definition \ref{Definition operator}, we write $\tilde{\mathcal L}$  equivalently as
\begin{align*}
	\tilde{\mathcal L}&=y^{ \alpha}\mbox{Tr }\left(QD^2_x\right)+2y^{ \alpha}\left(\tilde q,\nabla_xD_y\right)+\tilde \gamma y^{ \alpha}  \tilde L_y+y^{ \alpha-1}\left(d,\nabla_x\right)
\end{align*}
where  
\begin{align*}
\tilde L_y=D_{yy}+\frac{\tilde c/\tilde \gamma}{y}D_y-\frac{b/\tilde \gamma}{y^2}.
\end{align*}
  
By \eqref{tilde D gamma}, \eqref{tilde s gamma} and observing that, by the assumption on $\alpha_1$, $\alpha_2$, ${\beta_\alpha}+1=\frac{\alpha_1-\alpha_2+2}{2}>0$, the discriminant $\tilde D$ and the parameters $\tilde s_{1,2}$ of $\tilde L_y$ defined as in \eqref{defD}, \eqref{defs}  are related to  those of $L_y$ by
	\begin{align*}
		\tilde D&=\frac{D}{({\beta_\alpha}+1)^2},\qquad \tilde s_{1,2}=\frac{s_{1,2}}{{\beta_\alpha}+1}
	\end{align*}

 The generation properties, the domain description and the  maximal regularity  for the operator in  \eqref{operator alpha diff def} can be then deduced by  the same properties for the operator studied in the previous sections  when  $\alpha_1=\alpha_2$.\\

 We start by the case $b=0$ with Neumann and oblique boundary condition. We recall that the   Sobolev spaces $W^{2,p}_{v}\left(\alpha_1,\alpha_2,m\right)$ and  $W^{2,p}_{\mathcal N}\left(\alpha_1,\alpha_2,m\right)$, as well as their boundary conditions,  are introduced and analysed  in details in Sections \ref {Section neumann sobolev} and  \ref{sobolev oblique alpha diff section} to which we refer.

\begin{teo} \label{complete-oblique}
Let $v=(d,c)\in\R^{N+1}$ with $d=0$ if $c +\beta_\alpha\gamma=0$, and  	let $\alpha_1,\alpha_2\in\R$ such that $\alpha_2<2$,  $\alpha_2-\alpha_1<2$. If 
	$$\alpha_1^{-} <\frac{m+1}p<\frac{c}{\gamma}+1-\alpha_2$$   
	then the operator
	\begin{align*}
		\mathcal L=y^{\alpha_1}\mbox{Tr }\left(QD^2_x\right)+2y^{\frac{\alpha_1+\alpha_2}{2}}q\cdot \nabla_xD_y+\gamma y^{\alpha_2}D_{yy}+y^{\frac{\alpha_1+\alpha_2}{2}-1}\left(d,\nabla_x\right)+c y^{\alpha_2-1}D_y
	\end{align*}
	 endowed with domain
 \begin{align*}
	&W^{2,p}_{w}\left(\alpha_1,\alpha_2,m\right),\qquad w=(d,c+{\beta_\alpha}\gamma)
\end{align*}  
  generates a bounded analytic semigroup  in $L^p_m$ which has maximal regularity.  Moreover the set  $ \mathcal C_w$ defined in \eqref{defC oblique} is a core for $\mathcal L$.
\end{teo}
{\sc Proof.}  According to the discussion above  the isometry 
$$T_{0,{\beta_\alpha}}:L^p_{\tilde{m}}\to L^p_m,\qquad  {\beta_\alpha}=\frac{\alpha_1-\alpha_2}{2},\qquad \tilde m=\frac{m-{\beta_\alpha}}{{\beta_\alpha}+1}$$
transforms  $\mathcal L$ into
\begin{align*}
	\tilde{\mathcal L}=y^{\alpha}\mbox{Tr }\left(\tilde AD^2\right)+y^{\alpha-1}\left(\tilde v,\nabla\right)
\end{align*}
where $\tilde A$, $\tilde v$ are defined in \eqref{coeff tilde L eq} and $\tilde{\mathcal L}$ acts on $L^p_{\tilde m}$.
Observe now that the assumptions on the parameters translates into  $\alpha<2$ and $\alpha^- < \frac{\tilde m +1}{p} < \frac{\tilde c}{\tilde \gamma }+1-\alpha$. The statement on generation and maximal regularity is therefore a translation to  $\mathcal L$ and in $L^p_m$ of the results of Theorem \ref{Teo alpha Oblique} for $\tilde{\mathcal L}$ in $L^p_{\tilde m}$.
Concerning the domain, we have
$$D(\mathcal L)=T_{0,{\beta_\alpha}}\left( W^{2,p}_{\tilde v}\left (\alpha, \alpha,  \tilde m \right)\right)$$
which, by Proposition \ref{Isometry action sob}, coincides with $W^{2,p}_{ w}\left(\alpha_1,\alpha_2,m\right)$.\qed

\begin{os} As in the proof of Theorem \ref{Teo alpha Oblique}, when $c +\beta_\alpha\gamma\neq 0$  we can  transform $\mathcal L$  into a similar operator with $d=0$ and Neumann boundary condition. 
	Indeed,  we  use the the isometry  \eqref{Tran def} with $\omega=-\frac{d}{(c+\gamma\beta_\alpha)(\beta_\alpha+1)}$, namely
	\begin{align*}
		S_{\beta_\alpha,\omega}\, u(x,y)&:=u\left(x-\omega y^{\beta_\alpha+1},y\right),\quad (x,y)\in\R_+^{N+1}.
	\end{align*}
	Then by Corollary \ref{Isometry action 2 spe} and by Proposition \ref{Equi obli Neumm sob} one has 
\begin{align*}
	{S_{\beta_\alpha,\omega}}^{-1} \Big(\mathcal L\Big){S_{\beta_\alpha,\omega}}
	&=\tilde{\mathcal L},\qquad D(\tilde{\mathcal L})=W^{2,p}_{\mathcal N}\left (\alpha_1, \alpha_2,  m \right)
\end{align*}
	where
	\begin{align*}
		\tilde{\mathcal L}=y^{\alpha_1}\mbox{Tr }\left(\tilde QD^2_x\right)+2y^{\frac{\alpha_1+\alpha_2}2}\left(\tilde q,\nabla_xD_y\right)+\gamma y^{\alpha_2}D_{yy}+c y^{\alpha_2-1}D_y
	\end{align*}
and  \;$\tilde Q=Q-2q\otimes \frac{d}{c+\gamma\beta_\alpha}+\gamma\frac{d\otimes d}{(c+\gamma\beta_\alpha)^2}$,\; $\tilde q=q-\frac{\gamma}{c+\gamma\beta_\alpha}d$.
\end{os}
\medskip
 
	We point out, without stating them explicitly, that analogous results as those in  Corollaries \ref{Oblique cor1}, \ref{Oblique cor2} apply also to this case. \\
	
\medskip 
We finally add  the potential term  and study $\mathcal L$  under Dirichlet boundary conditions.

\begin{teo} \label{complete dirichlet alpha diff} Let  $\alpha_2<2$, $\alpha_2-\alpha_1<2$ and
	$$s_1+ \alpha_1^-<\frac{m+1}p<s_2+2-\alpha_2.$$   
	Then the operator
	\begin{align*}
		\mathcal L=y^{\alpha_1}\mbox{Tr }\left(QD^2_x\right)+2y^{\frac{\alpha_1+\alpha_2}{2}}q\cdot \nabla_xD_y+\gamma y^{\alpha_2}D_{yy}+y^{\frac{\alpha_1+\alpha_2}{2}-1}\left(d,\nabla_x\right)+c y^{\alpha_2-1}D_y-by^{\alpha_2-2}
	\end{align*} 
	generates a bounded analytic semigroup  in $L^p_m$ which has maximal regularity. Moreover,
	\begin{equation*}
		D(\mathcal L)
		=y^{-s_1}W^{2,p}_{w}\left(\alpha_1,\alpha_2,m-s_1p\right),\qquad  w=(d,c+\beta_\alpha\gamma)-2s_1\left(q,\gamma\right).
	\end{equation*}	
\end{teo}
{\sc Proof.}  As in the proof of the previous theorem,  the isometry 
$$T_{0,{\beta_\alpha}}:L^p_{\tilde{m}}\to L^p_m,\qquad  {\beta_\alpha}=\frac{\alpha_1-\alpha_2}{2},\qquad \tilde m=\frac{m-{\beta_\alpha}}{{\beta_\alpha}+1}$$
transforms  $\mathcal L$ into
\begin{align*}
		\tilde{\mathcal L}&=y^{\alpha}\mbox{Tr }\left(\tilde AD^2\right)+y^{\alpha-1}\left(\tilde v,\nabla\right)- by^{\alpha-2}\\[1ex]
		&=y^{ \alpha}\mbox{Tr }\left(QD^2_x\right)+2y^{ \alpha}\left(\tilde q,\nabla_xD_y\right)+\tilde \gamma y^{ \alpha}  \tilde L_y+y^{ \alpha-1}\left(d,\nabla_x\right)
	\end{align*}
where $\tilde A$, $\tilde v$ are defined in \eqref{coeff tilde L eq} and $\tilde{\mathcal L}$ acts on $L^p_{\tilde m}$. Moreover   the parameters $\tilde s_{1,2}$ of $\tilde L_y$ satisfies
\begin{align}\label{parameter last section}
\tilde s_{1,2}=\frac{s_{1,2}}{{\beta_\alpha}+1}.
\end{align}
Observe now that the hypotheses on the parameters translates into  $\alpha<2$ and
$$\tilde{s_1}+ \alpha_1^-<\frac{\tilde m+1}p<\tilde{s_2}+2-\alpha_2.$$
 The statement on generation and maximal regularity is therefore a translation to  $\mathcal L$ and in $L^p_m$ of the results of Theorem \ref{complete} for $\tilde{\mathcal L}$ in $L^p_{\tilde m}$.
Concerning the domain we have
\begin{align*}
	D(\mathcal L)&=T_{0,{\beta_\alpha}}\left(D(\tilde{\mathcal L})\right)=T_{0,{\beta_\alpha}}\Big(y^{-\tilde{s_1}} W^{2,p}_{\tilde w}\left (\alpha, \alpha,  \tilde m-\tilde{s_1}p \right)\Big)
\end{align*}
where  
\begin{align*}
	\tilde w=\tilde v-2\tilde{s_1}\left(\tilde q,\tilde\gamma\right)=\Big(d,(c+\beta_\alpha\gamma)(\beta_\alpha+1)\Big)-2s_1\big(q,\gamma(\beta_\alpha+1)\big).
\end{align*}
Then recalling \eqref{parameter last section} and using property (ii)-1 of  Proposition \ref{Isometry action der} and (ii) of Proposition \ref{Isometry action sob} we get
\begin{align*}
	D(\mathcal L)&=y^{-s_1}T_{0,{\beta_\alpha}}\Big( W^{2,p}_{\tilde w}\left (\alpha, \alpha,  \tilde m-\tilde{s_1}p \right)\Big)=y^{-s_1}W^{2,p}_{ w}\left(\alpha_1,\alpha_2,m-s_1p\right).
\end{align*}\qed

As in Corollary \ref{cor1}  we can characterized $D(\mathcal L)$ through a Dirichlet boundary condition.

\begin{cor} \label{corollary dirichelt alpha diffe equiv charct}
 Let  $\alpha_2<2$, $\alpha_2-\alpha_1<2$ such that
$$s_1+ \alpha_1^-<\frac{m+1}p<s_2+2-\alpha_2.$$   
Let $\mathcal A$ be the operator
\begin{align*}
	\mathcal A:&=\gamma y^{\alpha_2}D_{yy}+y^{\frac{\alpha_1+\alpha_2}2-1} (d-2s_1q) \cdot \nabla_x + c y^{\alpha_2-1}D_y-by^{\alpha_2-2}.
\end{align*}
	\begin{itemize}
		\item[(i)] If $D>0$ then 
		\begin{align*} 
			D(\mathcal L)=&\left\{u \in  W^{2,p}_{loc}(\R^{N+1}_+): u,\ y^{\alpha_1}\Delta_xu,\ \mathcal Au \in L^p_m\text{\;\;and\;\;}\lim_{y\to 0}y^{s_2} u=0\right\}
		\end{align*}
		\item[(ii)] If $D=0$ then $s_1=s_2$ and
		\begin{align*} 
			D(\mathcal L)=&\left\{u \in  W^{2,p}_{loc}(\R^{N+1}_+): u,\ y^{\alpha_1}\Delta_xu,\ \mathcal Au\in L^p_m\text{\;\;and\;\;}\lim_{y\to 0}y^{s_2} u\in\C\right\}
		\end{align*}
	\end{itemize} 
	In both cases the graph   norm and $$\|u\|_{L^p_m}+\|y^{\alpha}\Delta_x u\|_{L^p_m}+\|\mathcal Au\|_{L^p_m}$$ are equivalent on $D(\mathcal L)$. 
\end{cor}
{\sc Proof.} Let us prove claim (i). By Theorem \ref{complete dirichlet alpha diff} we have 
\begin{equation*}
	D(\mathcal L)
	=T_{-s_1,0}\left(W^{2,p}_{w}\left(\alpha_1,\alpha_2,m-s_1p\right)\right)
\end{equation*}	
where
\begin{equation*}
 w=(d,c+\beta_\alpha\gamma)-2s_1\left(q,\gamma\right):=\left(w_x,w_{N+1}\right).
\end{equation*}	
Then we proceed as in the proof of Corollary \ref{cor1}: we apply Proposition Proposition \ref{Equiv Wv alpha diff dir}  with  $\mu=\frac{c}{\gamma}-2s_1=1+2\sqrt{D}> 1$ obtaining
	\begin{align*} 
	W^{2,p}_{w}(\alpha_1, \alpha_2, m-s_1p)=&\left\{v \in  W^{2,p}_{loc}(\R^{N+1}_+): v,\ y^{\alpha_1}\Delta_xv,\   	\tilde{\mathcal A}v \in L^p_{m-s_1p}\text{\;\;and\;\;}\lim_{y\to 0}y^{\mu-1} v=0\right\}
\end{align*}
where
\begin{align*}
	\tilde{\mathcal A}:&=y^{\alpha_2}D_{yy}+y^{\frac{\alpha_1+\alpha_2}2-1} \tilde d \cdot \nabla_x + \mu y^{\alpha_2-1}D_y,\qquad \tilde d=\frac{\mu+\beta_\alpha}{w_{N+1}}w_x=\frac{d-2s_1q}{\gamma}.
\end{align*}
The required claim then follows from the previous equalities after noticing that, by Proposition \ref{Isometry action} and Corollary \ref{Isometry action2}, 
\begin{align*}
	&T_{-s_1,0}^{-1}\left( y^{\alpha_1} \Delta_x\right) T_{-s_1,0}=y^{\alpha_1}\Delta_x,\qquad  T_{-s_1,0}^{-1}\left(\mathcal A \right)T_{-s_1,0}=\gamma \tilde{\mathcal A}.
\end{align*} 
and that  for any $u\in y^{-s_1}W^{2,p}_{w}\left(\alpha_1,\alpha_2,m-s_1p\right)$, setting $v=y^{s_1}u$, one has, recalling \eqref{defs},
\begin{align*}
	y^{\mu -1}v=y^{\mu-1+s_1}u=y^{s_2 }u.
\end{align*} 
Claim (ii) follows similarly.
\qed

As in Corollary \ref{CorBdalpha}, in some range of parameters one has an   improvement in the elliptic regularity  of  the operator.
\begin{cor} 	Let  $\alpha_2<2$, $\alpha_2-\alpha_1<2$ and
	$$s_1+ \alpha_1^-<\frac{m+1}p<s_2+2-\alpha_2.$$ 
	\begin{itemize}
		\item[(i)] If both the condition  $s_1+ \alpha_1^-<\frac{m+1}p<s_2+2-\alpha_2$ and  $\frac{m+1}p>s_1+1-\frac{\alpha_1+\alpha_2}2$ hold then $D(\mathcal L)
		=y^{-s_1}W^{2,p}_{\mathcal N}\left(\alpha_1,\alpha_2,m-s_1p\right)$ and 
		\begin{equation*} 
			\|y^{\alpha_1} D_{x_i x_j} u\|_{L^p_{m}} +\| y^{\alpha_2} \gamma L_yu\|_{L^p_{m}}+\|y^{\frac{\alpha_1+\alpha_2}2} D_{x_iy} u\|_{L^p_{m}}+\|y^{\frac{\alpha_1+\alpha_2}2-1}  \nabla_x u\|_{L^p_{m}}\leq C\| \mathcal Lu\|_{L^p_{m}}
		\end{equation*}
		where $ y^{\alpha_2} \gamma L_y=y^{\alpha_2} \gamma D_{yy} +cy^{\alpha_2-1}D_y-by^{\alpha_2-2}$.
		\item[(ii)] If $s_1+2-\alpha_2<\frac{m+1}p<s_2+2-\alpha_2$ then $D(\mathcal L)=W^{2,p}_{\mathcal R}(\alpha_1,\alpha_2,m)$.
	\end{itemize}
\end{cor}
{\sc{Proof.}} Identical to the proof of Corollary \ref{Cor-Rellich}. For claim (ii) we also observe that  $s_1+ 2-\alpha_2>s_1+\alpha_1^- $, since  $\alpha_2<2$, $\alpha_2-\alpha_1<2$. \qed

\begin{os}
When $b=0$, we remark, without stating explicitly, that   the  results of Remark \ref{oss uniquenss DirNeu} and Corollary \ref{CorBdalpha}  apply also in this case . In particular if  $c>\gamma $ (so that $s_1=0$, $s_2=\frac c\gamma -1$), the operator  of Theorems \ref{complete-oblique} and \ref{complete dirichlet alpha diff} coincide, showing that in this case the Dirichlet and the oblique derivative boundary conditions are equivalent. On the other hand if  $c<\gamma $ (so that $s_1=\frac c\gamma -1 \neq 0$, $s_2=0$) we can construct a realization of $\mathcal L$  different from that of Theorem  \ref{complete-oblique}.
\end{os}

\bibliography{../../../TexBibliografiaUnica/References}

\end{document}